\documentclass[12pt]{amsart}

\makeatletter

\providecommand{\LyX}{L\kern-.1667em\lower.25em\hbox{Y}\kern-.125emX\@}

 \theoremstyle{plain}    
 \newtheorem{thm}{Theorem}[section]
 \numberwithin{equation}{section} %% Comment out for sequentially-numbered
 \numberwithin{figure}{section} %% Comment out for sequentially-numbered
 \theoremstyle{plain}    
 \newtheorem{cor}[thm]{Corollary} %%Delete [thm] to re-start numbering
 \theoremstyle{plain}    
 \newtheorem{lem}[thm]{Lemma} %%Delete [thm] to re-start numbering
 \theoremstyle{plain}    
 \newtheorem{prop}[thm]{Proposition} %%Delete [thm] to re-start numbering
 \theoremstyle{definition}
 \newtheorem{defn}[thm]{Definition}
 \theoremstyle{remark}
 \newtheorem{rem}[thm]{Remark}
 \theoremstyle{remark}    
 \newtheorem{claim}[thm]{Claim}

\newcommand{\Tor}{\operatorname{Tor}}
\newcommand{\Image}{\operatorname{im}}
\newcommand{\Span}{\operatorname{span}}
\newcommand{\Hom}{\operatorname{Hom}}
\newcommand{\Tr}{\operatorname{Tr}}

\def\n{{\noindent}}
\def\g{{\gamma}}

\makeatother
\begin{document}

\baselineskip=0.5cm 
 \thispagestyle{empty}
\begin{center}
{\large\bf $L^{2}$-HOMOLOGY FOR VON NEUMANN ALGEBRAS}
\end{center}

\vspace{0.5cm}

\begin{center} Alain CONNES \footnote{Coll\`ege de France, 3 rue
d'Ulm, 75  005 Paris, \\ I.H.E.S. and Vanderbilt University,
connes$@$ihes.fr}
and 
 Dimitri SHLYAKHTENKO
\footnote{Department of Mathematics University of California, Los Angeles, CA 90095, USA\\ shlyakht@math.ucla.edu. 
Reserach supported by the NSF and the Sloan Foundation}
\end{center} \vspace{0,25cm}

\section{Introduction.}

The aim of this paper is to introduce a notion of $L^{2}$-homology
in the context of von Neumann algebras. Finding a suitable (co)homology
theory for von Neumann algebras has been a dream for several generations
(see \cite{kad-ringI,kad-ringII,kad-ringIII,sinclair-smith:cohomology}
and references therein). One's hope is to have a powerful invariant
to distinguish von Neumann algebras. Unfortunately, little positive
is known about the Kadison-Ringrose cohomology $H_{b}^{*}(M,M)$,
except that it vanishes in many cases. Furthermore, there does not
seem to be a good connection between the bounded cohomology theory
of a group and of the bounded cohomology of its von Neumann algebra.

\n Our interest in developing an $L^{2}$-cohomology theory was revived
by the introduction of $L^{2}$-cohomology invariants in the field
of ergodic equivalence relations in the paper of Gaboriau
\cite{gaboriau:ell2}. His results in particular imply that $L^{2}$-Betti
numbers $\beta _{i}^{(2)}(\Gamma )$ of a discrete group are the same
for measure-equivalent groups (i.e., for groups that can generate
isomorphic ergodic measure-preserving equivalence relations). Parallels
between the {}``worlds'' of von Neumann algebras and measurable
equivalence relations have been noted for a long time (starting with
the parallel between the work of Murray and von Neumann \cite{vN} and that of H. Dye
\cite{Dye}).
% to notions of amenability in both of these
% worlds \cite{connes:injective,connes:feldman:weiss:amenable}).
 Thus
there is hope that an invariant of a group that {}``survives'' measure
equivalence will survive also {}``von Neumann algebra equivalence'',
i.e., will be an invariant of the von Neumann algebra of the group.

\n  The original motivation for our construction comes from 
the well understood analogy between the theory of II$_1$-factors
and that of discrete groups, based on the theory of correspondences 
\cite{connes:correspondences,connes:ncgeom}. This analogy has been remarkably
efficient to transpose analytic properties such as ``amenability" or ``property $T$"
from the group context to the factor context \cite{connesT} \cite{cj}
and more recently in the breakthrough work of Popa \cite{popa} \cite{bbk}. 
We use the theory of correspondences 
 together with
the algebraic
description of $L^{2}$-Betti numbers given by Luck. His definition
involves the computation of the algebraic group homology with coefficients
in the group von Neumann algebra. Following the guiding
idea that the category of bimodules over a von Neumann algebra is
the analogue of the category of modules over
a group, we are led to the following algebraic definition of $L^{2}$-homology
of a von Neumann algebra $M$:\[
H_{k}^{(2)}(M)=H_{k}(M;M\bar{\otimes }M^{o}).\]
Here $H_{k}$ stands for the algebraic Hochschild homology of $M$.
One is thus led to consider the $L^{2}$-Betti numbers,\[
\beta _{k}^{(2)}(M)=\dim _{M\bar{\otimes }M^{o}}H_{k}^{(2)}(M)\]
(see Section \ref{sub:Fromgroupstovnalg} for more motivation behind
this definition).

\n The $L^2$ Betti numbers 
that we associate to a II$_1$ factor $M$ enjoy the following scaling property.  
If $p\in M$ is a projection of trace $\lambda$, then $$
\beta_k^{(2)}(pMp) = \frac {1 }{\lambda ^2} \beta_k^{(2)} (M).$$  
This behavior is  a consequence of considering the ``double'' $M\otimes M^o$ of $M$.
It is quite different from the behavior of $L^2$-Betti numbers associated to 
equivalence relations \cite{gaboriau:ell2}, where the factor $\lambda^2$ is replaced by
$\lambda$.
Note that the behavior of our invariants is consistent with the formula for the 
``number of generators'' of compressions of free group factors,
where the same factor $\lambda^2$ appears.  Indeed,  for $p\in L(\mathbb{F}_t)
= L(\mathbb{F}_{(t-1)+1})$
a projection of trace $\lambda$, one has  by the results of 
Voiculescu, Dykema and Radulescu (see e.g. \cite{dvv:book}) that $$pL(\mathbb{F}_t) p
\cong L(\mathbb{F}_{ ((t-1)\lambda ^{-2}) + 1}.$$
In this respect, our theory seems to be quite disjoint from a kind of relative 
theory of $L^2$-invariants for factors with $HT$-Cartan subalgebras discovered by 
Popa \cite{popa}. Indeed, the two theories have different scaling properties ($\lambda^{-2}$
versus $\lambda^{-1}$) under compression.

\medskip
\n We have completely unexpectedly found
a connection between our definition, involving algebraic homology,
and the theory of free entropy dimension developed by Voiculescu and
his followers (\cite{dvv:entropy1,dvv:entropy2,dvv:entropy3,dvv:entropy5},
see also \cite{dvv:entropysurvey} and references therein). Such a
connection between $L^{2}$-Betti numbers and free entropy dimension
parallels nicely the connections between free entropy dimension and
cost for equivalence relations 
\cite{gaboriau:cost,shlyakht:cost,shlyakht:cost:micro}. 

\n Our definition seems to be robust with respect to modifications; for example,
the connections with free probability theory do not seem to be accessible 
for certain variations of our definition.

\n This unexpected connection allows us to rely on several deep results
in free probability theory, including the inequality between the microstates
and microstates-free free entropy proved by recent fundamental work
of Biane, Capitaine and Guionnet \cite{guionnet-biane-capitaine:largedeviations}.
Drawing on these results, we obtain strong evidence that the first
$L^{2}$-Betti number associated to the von Neumann algebra of a free
group $\mathbb{F}_{n}$ on $n$ generators is bounded below by $n-1$
($=\beta _{1}^{(2)}(\mathbb{F}_{n})$). Moreover, our inability to
prove a similar inequality between free entropy dimension and $L^{2}$-Betti
numbers obtained by varying our definition can be taken as evidence
that the present definition is the correct one. Even more encouraging
are the facts that the $L^{2}$-Betti numbers  can be controlled
from above at least in some cases (such as abelian von Neumann algebras,
or von Neumann algebras with a diffuse center).

\n Our results also have several implications for the questions involving
computation of free entropy dimension. In particular, we show that
the modified free entropy dimension $\delta _{0}(\Gamma )$ of any
discrete group $\Gamma $ is bounded from above by $\beta _{1}^{(2)}(\Gamma )-\beta
_{0}^{(2)}(\Gamma )+1$;
a similar estimate holds for the various versions of the non-microstates
free entropy dimension. If $\Gamma $ has property $T$ and if the
von Neumann algebra $L(\Gamma )$ is diffuse, then $\delta _{0}(\Gamma )\leq 1$
(and $=1$ if $L(\Gamma )$ can be embedded into the ultrapower of
the free group factor, as is the case when $\Gamma $ is residually
finite). This generalizes a result of Voiculescu \cite{dvv:entropySLnZ}
for the case that $\Gamma =SL(n,\mathbb{Z})$, $n\geq 3$.

\subsection{Notation.}

Whenever possible, we shall use the following notation: the letter
$A$ will stand for a tracial $\ast$-algebra with a fixed positive faithful
trace $\tau $. The letter $M$ will stand for the von Neumann algebra
generated by $A$ in the GNS representation associated to $\tau $.
We write $L^{2}(M)$ for the associated representation space. 

\n For a group $\Gamma $, we denote by $\mathbb{C}\Gamma $ its (algebraic)
group algebra, and by $L(\Gamma )$ its group von Neumann algebra.

\n The tensor sign $\otimes $ will always refer to the algebraic tensor
product. Tensor products that involve completions (such as the von
Neumann algebra tensor product or the Hilbert space tensor product)
will always be denoted by $\bar{\otimes }$.

\n We denote by $HS(L^{2}(M))$ the ideal of Hilbert-Schmidt operators
on $L^{2}(M)$, and by $FR(L^{2}(M))$ the ideal of finite-rank operators
on $L^{2}(M)$.

\n It will be useful to identify $HS(L^{2}(M))$ with $L^{2}(M)\bar{\otimes }L^{2}(M^{o})$
in the following way. Let $J:L^{2}(M)\to L^{2}(M)$ be the Tomita
conjugation. Let $P_{1}\in L^{2}(M)$ be the projection onto the cyclic
vector. Then consider the map $\Psi :M\otimes M^{o}\to FR(L^{2}(M))$
given by\begin{equation}
\Psi (x\otimes y^{o})(\xi )=x\  \tau (y\,\xi )\qquad x,y\in M,\  \xi \in
L^{2}(M).\label{eq:defofPsi}\end{equation}
The map $\Psi $ extends to an isometric isomorphism of $L^{2}(M\bar{\otimes }M^{o})$
with the Hilbert space $HS(L^{2}(M))$ of Hilbert-Schmidt operators
in $L^{2}(M)$. 

\n Note that $B(L^{2}(M))$ admits four commuting actions of $M$: \begin{eqnarray*}
T & \mapsto  & aT\\
T & \mapsto  & Ta\\
T & \mapsto  & Ja^{*}JT\\
T & \mapsto  & TJa^{*}J.
\end{eqnarray*}
Here $a\in M$ and $T\in B(L^{2}(M))$. These four actions are intertwined
by the map $\Psi $ with the four natural actions of $M$ on $L^{2}(M\bar{\otimes
}M^{o})$,
listed in the corresponding order:\begin{eqnarray}\label{bim}
x\otimes y^{o} & \mapsto  & a\,x\otimes y^{o}\\
x\otimes y^{o} & \mapsto  & x\otimes (y\,a)^{o} \nonumber\\
x\otimes y^{o} & \mapsto  & x\,a\otimes y^{o}\nonumber\\
x\otimes y^{o} & \mapsto  & x\otimes (a\,y)^{o},\nonumber
\end{eqnarray}
where $a\in M$ and $x\otimes y^{o}\in M\otimes M^{o}$.

\section{Definition of $L^{2}$-Betti numbers and $L^{2}$-homology.}

\subsection{Review of $L^{2}$-homology for groups.}

\subsubsection{Co-compact actions.}

The study of $L^{2}$-cohomology has been initiated by Atiyah in \cite{atiyah-L2},
who considered $L^{2}$ deRham cohomology of a connected manifold
$X$ endowed with a cocompact proper free action of a discrete group $\Gamma $.
The $k$-th cohomology group is defined as the quotient of the space
of closed $L^{2}$-integrable $k$-forms by the closure of the space
of exact $L^{2}$ $k$-forms. In other words, one considers \emph{reduced}
$L^{2}$-cohomology (reduced means that one takes a quotient by the
closure of the image of the boundary operator).

\n The cohomology groups in question can also be identified with the
Hilbert spaces of $L^{2}$ harmonic forms, and are in a natural way
modules over the group $\Gamma $. Furthermore, because of the properness
and cocompactness assumptions, these spaces are actually  modules over
the group von Neumann algebra $L(\Gamma )$. Atiyah considered the
Murray-von Neumann dimensions of these cohomology groups, and called
them $L^{2}$-Betti numbers of the action. Remarkably, it turns out
that these numbers are homotopy invariants; in particular, if the
manifold $X$ is contractible, then the numbers depend only on the
group $\Gamma $, and are called the $L^{2}$-Betti numbers of the
group, $\beta _{i}^{(2)}(\Gamma )$. This definition extends to cover 
cocompact proper free actions of  $\Gamma $ on a contractible 
simplicial complex, and gives an equivalent definition of the 
Betti numbers using $L^{2}$ singular
homology. For example, if $\Gamma =\mathbb{F}_{n}$,
the free group on $n$ generators, then $\beta _{i}^{(2)}(\Gamma )=0$
for $i\neq 1$, and $\beta _{1}^{(2)}(\Gamma )=n-1$.

\subsubsection{Cheeger-Gromov's approach to non-cocompact actions.}

\n Not every group can act on a contractible space in a proper free and co-compact
way. For this reason, if one wants to obtain an $L^{2}$-homology
theory for groups, one is forced to consider non-cocompact actions.
This leads to difficulties of a technical nature. The problem is that
because the spaces of $L^{2}$-(co)chains are now infinitely-dimensional
modules over $L(\Gamma )$, one needs to be much more careful in taking
the closure of the image of the boundary operator, and one may end
up with having to consider dimensions of actions of $L(\Gamma )$
on quotients of Hilbert spaces by not necessarily closed subspaces.

\n This can be overcome in one of two ways. The first, which is the original
approach of Cheeger and Gromov \cite{cheeger-gromov:l2}, is to {}``approximate''
the $L^{2}$-homology of a non-cocompact action on a contractible
manifold $X$ by realizing the simplicial complex $(\Delta _{*},\partial _{*})$
underlying $X$ as an inductive limit of sub-complexes $(\Delta _{*}^{(k)},\partial
_{*}^{(k)})$,
which correspond to co-compact (free) actions. Let $C_{*}^{(k)}(X)$
the the space of $L^{2}$-chains for $\Delta _{*}^{(k)}$ (i.e., the
 Hilbert space with orthonormal basis given by simplices in $\Delta _{*}^{(k)}$).
Denote by $i_{k,l}$ the map from $C_{*}^{(k)}(X)$ to $C_{*}^{(l)}(X)$.
Then the $n$-th $L^{2}$-Betti number is defined as\begin{equation}
\beta _{n}^{(2)}=\sup _{k\to \infty }\inf _{l>k}\dim _{L(\Gamma )}\frac{i_{k,l}(\ker
\partial _{n}^{(k)})}{\overline{i_{k,l}(\ker \partial _{n}^{(k)})\cap \Image
\partial _{n+1}^{(l)}}}.\label{eq:l2bettinoncocompactChGr}\end{equation}

\subsubsection{Luck's approach to non-cocompact actions.}

The second way, developed by Luck \cite{luck:foundations1}, is to
extend the notion of Murray-von Neumann dimension to algebraic modules
over type II$_{1}$ von Neumann algebras. This way one can, for example,
assign a dimension to the quotient of an $L(\Gamma )$ module by a
non-closed submodule. Luck shows that such an extension indeed exists.
In the case that $V$ is a finitely-generated module over $L(\Gamma )$,
its dimension is just the supremum of the dimensions of \emph{normal}
$L(\Gamma )$ modules (i.e. finite projective modules) that can be embedded into $V$. 

\n Returning briefly to the co-compact case,  one can now consider
\emph{non-reduced} simplicial $L^{2}$ homology of a space $X$,
defined as the quotient of the kernel of the boundary operator $\partial _{n}$
by the (non-closed) image of $\partial _{n+1}$. This results in $L(\Gamma )$
modules, which are not normal. However, due to the behavior of Luck's
extension of the Murray-von Neumann dimension, it turns out that the
dimensions of these modules are the same as the Murray-von Neumann
dimensions of the reduced homology groups (i.e., the $L^{2}$-Betti
numbers). 

\n In the case that $X$ is connected and contractible, its ordinary
homology vanishes. This means that if we denote by $C_{k}^{(f)}$
the vector space with basis given by the $k$-chains on $X$, then
the sequence $(C_{k}^{(f)},\partial _{k})$ is exact. Each $C_{k}^{(f)}$
is a flat module over $\Gamma $, and $C_{0}^{(f)}\cong \mathbb{C}$,
since we assume that $X$ is connected. Thus $(C_{k}^{(f)},\partial _{k})$
forms a resolution of the trivial $\Gamma $-module $C_{0}^{(f)}$.
Furthermore, the space of $L^{2}$ chains on $X$ is very roughly
$L(\Gamma )\otimes _{\Gamma }C_{k}^{(f)}$.%
\footnote{More precisely, the space of $L^{2}$-chains is $\ell ^{2}(\Gamma )\otimes
_{\Gamma }C_{k}^{(f)}=\ell ^{2}(\Gamma )\otimes _{L(\Gamma )}L(\Gamma )\otimes
_{\Gamma }C_{k}^{(f)}$.
However, since by \cite{luck:foundations1}, the functor $\ell ^{2}\otimes _{L(\Gamma
)}\cdot $
in the category of finitely-generated $L(\Gamma )$-modules is flat,
the nuance between $\ell ^{2}(\Gamma )\otimes _{\Gamma }C_{k}^{(f)}$
and $L(\Gamma )\otimes _{\Gamma }C_{k}^{(f)}$ is irrelevant in the
foregoing.%
} Thus the non-reduced $L^{2}$-homology of $X$ is nothing but the
(algebraic) homology group\[
H_{*}(\Gamma ;L(\Gamma ))\]
of the group $\Gamma $ with coefficients in $L(\Gamma )$ viewed as
a right $\Gamma $-module. Since
$L(\Gamma )$ is both a right $\Gamma $-module, and a left $L(\Gamma )$-module,
the $H_{*}(\Gamma ;L(\Gamma ))$ are left $L(\Gamma )$-modules. Thus
one can consider $\dim _{L(\Gamma )}H_{*}(\Gamma ;L(\Gamma ))$. However,
as we pointed out, the dimensions of these non-reduced homology groups
are precisely the $L^{2}$-Betti numbers of the group.

\n Thus even in the case that $\Gamma $ admits no proper free co-compact action
on a contractible space, one can consider the numbers\begin{equation}
\beta _{*}^{(2)}(\Gamma )=\dim _{L(\Gamma )}H_{*}(\Gamma ;L(\Gamma
)),\label{eq:l2bettiLuck}\end{equation}
where $H_{*}$ stands for the algebraic group homology. Luck shows
that a certain continuity property of his dimension leads to the formula\[
\dim _{L(\Gamma )}H_{*}(\Gamma ;L(\Gamma ))=\sup _{k}\inf _{l>k}\dim _{L(\Gamma
)}i_{k,l}(H_{*}(C_{*}^{(k)}))\]
in the notation of the previous section and of equation
(\ref{eq:l2bettinoncocompactChGr}).
Combined with (\ref{eq:l2bettinoncocompactChGr}), this shows that
this definition produces the same numbers as the one in equation
(\ref{eq:l2bettinoncocompactChGr}).

\subsection{From groups to von Neumann algebras.\label{sub:Fromgroupstovnalg}}

\subsubsection{Correspondences.}

Let us recall the basis of the analogy between discrete groups
and II$_1$-factors. A correspondence from 
II$_1$-factors is simply a Hilbert space endowed with a structure
of (normal) bimodule. Correspondences on $M$ (i.e. $M$-$M$ Hilbert bimodules)
play the role of unitary 
representations, and the dictionary begins as follows,   
$$
\begin{matrix}
\hbox{Discrete Group $\Gamma$} &\hbox{II$_1$-Factor $M$} \cr
\cr
\hbox{Unitary Representation} &\hbox{$M$-$M$ Hilbert bimodule} \cr
\cr
\hbox{Trivial Representation} &\hbox{$L^2(M)$} \cr
\cr
\hbox{Regular Representation} &\hbox{Coarse Correspondence } \cr
\cr
\hbox{Amenability} &\hbox{$L^2(M)\subset_{\textrm{weakly}}L^{2}(M)\bar{\otimes }L^{2}(M^{o})$} \cr
\cr
\hbox{Property T} &\hbox{$L^2(M)$ isolated}
\cr
\end{matrix}
$$

\medskip
\n where the ``coarse correspondence" is given by the 
bimodule $L^{2}(M)\bar{\otimes }L^{2}(M^{o})$
 with bimodule structure given by the first
two lines of (\ref{bim}).

\subsubsection{Some homological algebra.}

In order to find a suitable notion of $L^{2}$ Betti numbers of a
tracial algebra $(M,\tau )$, we chose as our point of departure equation
(\ref{eq:l2bettiLuck}). It is perhaps best to rewrite the definition
of group homology in the language of homological algebra see e.g.
\cite{cartan-eilenberg:homalgebra}:\[
H_{*}(\Gamma ;L(\Gamma ))=\Tor _{*}^{Mod(\Gamma )}(\star ;L(\Gamma )).\]

\n This equation involves three objects: $Mod(\Gamma )$, $\star $ and
$L(\Gamma )$. The role of $Mod(\Gamma )$ is to prescribe a suitable
category; in this case, it is the category of left modules over $\Gamma $.
The module $\star $ stands for the trivial left module. To compute
the $\Tor $ functor, one must first choose a resolution of $\star $
by flat modules\begin{equation}
\star \leftarrow C_{0}\leftarrow C_{1}\leftarrow \cdots
.\label{eq:resolution}\end{equation}
This means that each $C_{i}$ is flat over $\Gamma $ (the definition
of a flat module is not really needed here; free modules and, more
generally, projective modules, are flat), and the sequence (\ref{eq:resolution})
is exact. 

\n Now $L(\Gamma )$  plays two roles. Its first role is that of a right
module over $\Gamma $. To finish the computation of $\Tor $, one
applies the functor $L(\Gamma )\otimes _{\Gamma }$ to the exact sequence
(\ref{eq:resolution}):\begin{equation}
L(\Gamma ) \otimes _{\Gamma }\star\leftarrow L(\Gamma
)\otimes _{\Gamma }C_{0}\leftarrow L(\Gamma
)\otimes _{\Gamma }C_{1}\leftarrow \cdots
.\label{eq:resolutiontensored}\end{equation}
This new sequence (\ref{eq:resolutiontensored}) may well fail to
be exact. Note, however, that it is a differential complex: the composition
of any two consecutive arrows is zero (since the original sequences
(\ref{eq:resolution}) had this property). The value of $\Tor $ is
precisely the homology of (\ref{eq:resolutiontensored}), i.e., $\Tor _{k}$
is the kernel of the $k$-th map, divided by the image of the $k+1$-st.
It is crucial that the value of the $\Tor $ functor is well-defined
and is independent of the choice of the resolution (\ref{eq:resolution}).

\n Finally, $L(\Gamma )$ plays its final role, which is that of not
just a right $\Gamma $-module, but also of a left $L(\Gamma )$-module.
It is crucial that the left action of $L(\Gamma )$ commutes
with the right action of $\Gamma $. Because of this, 
(\ref{eq:resolutiontensored}) is is a differential complex of 
left $L(\Gamma )$-modules. This
makes it possible to take the dimension
\begin{equation}
\beta _{*}^{(2)}(\Gamma )=\dim _{L(\Gamma )}\Tor _{*}^{Mod(\Gamma )}(\star ;L(\Gamma )),\label{eq:l2bettiLuck1}\end{equation}
 as the definition of the Betti
numbers (\ref{eq:l2bettiLuck}).

\subsubsection{Analogs for tracial $\ast$-algebras.}

Let now $(A,\tau )$ be a tracial $\ast$-algebra. We would like to make sense of
\[
\beta _{*}^{(2)}(A,\tau )=\dim _{X}\Tor _{*}^{Y}(Z,W),\]
where: $X$ is the analog of the group von Neumann algebra $L(\Gamma )$; $Y$ is
the analog of the category of right $\Gamma $-modules; $Z$ is the
analog of the trivial module; and $W$ is the analog of the $\Gamma ,L(\Gamma
)$-bimodule
$L(\Gamma )$.

\n Let $M=W^{*}(A)$ in the GNS representation associated to $\tau $.

\n It is fortunate that the theory of correspondences furnishes perfect
analogs for all of these objects. The analog of $Y$, i.e., 
of the category of left $\Gamma $-modules, is the category $Y$ of $A,A$-bimodules,
or, better, of left modules over the algebraic tensor product $A\otimes A^{o}$
of $A$ by its opposite algebra $A^{o}$. 

\n In particular the trivial $\Gamma $-module $Z$ is $A$ viewed as a bimodule over $A$
and since we want to view it as a {\it left} module over $A\otimes A^{o}$
we use
\begin{equation}\label{tri}
  (m \otimes n^{o})\cdot \,a := m\, a \, n \,, \quad \forall m\,, n\,, a \in A
\end{equation}

\n We next look for $X$, the analog of the group von Neumann algebra $L(\Gamma )$.
Note that $L(\Gamma )$ is the von Neumann algebra generated by $\Gamma $
in the left regular representation. The analog of the left regular
representation of $\Gamma $ on $\ell ^{2}(\Gamma )$ is the coarse
correspondence, i.e.
the representation of $A\otimes A^{o}$ on $L^{2}(M,\tau
)\bar{\otimes }L^{2}(M^{o},\tau )$ given by 
\[
(m\otimes n^{o})\cdot(\sum a_{i}\otimes b_{i}^{o}) =\sum m\, a_{i}\, \otimes \,(b_{i}\,n)^{o},\qquad m\,,n\in
A\,.\]
i.e. the first two lines of (\ref{bim}).
Hence $X=W^{*}(M\otimes M^{o})$ in this representation, i.e., $X=M\bar{\otimes }M^{o}$
(von Neumann tensor product). Finally, $W=M\bar{\otimes }M^{o}$,
but having two structures. One is that of a right $A\otimes A^{o}$-module, with
the action
\begin{equation}\label{bim1}
(\sum a_{i}\otimes b_{i}^{o})\cdot (m\otimes n^{o})=\sum a_{i}\, m\otimes (n\,b_{i})^{o},\qquad m\,,n\in
A\,.\end{equation}
The other is that of an $M\bar{\otimes }M^{o}$-left module, given by 
left multiplication in $M\bar{\otimes }M^{o}$, i.e. \[
(m\otimes n^{o})\cdot (\sum a_{i}\otimes b_{i}^{o})=\sum m\,a_{i}\otimes (b_{i}\,n)^{o}\,,\qquad
m,n\in M.\]
Note that the actions of $M\bar{\otimes }M^{o}$ and $A\otimes A^{o}$
on $W=M\bar{\otimes }M^{o}$ commute. 

\n Thus we are led by our analogy to consider
\begin{eqnarray*}
H_{k}^{(2)}(A,\tau ) & = & \Tor _{k}^{A\otimes A^{o}}(A;M\bar{\otimes }M^{o}),\\
\beta _{k}^{(2)}(A,\tau ) & = & \dim _{M\bar{\otimes }M^{o}}H_{k}^{(2)}(A,\tau ).
\end{eqnarray*}
It turns out that in the category of bimodules over an algebra, $\Tor _{k}^{A\otimes
A^{o}}(A;M\bar{\otimes }M^{o})$
is exactly the $k$-th Hochschild homology of $A$ with coefficients
in the bimodule $M\bar{\otimes }M^{o}$ of (\ref{bim1}). 

\begin{defn}
Let $(A,\tau )$ be a tracial $\ast$-algebra. Let $M=W^{*}(A)$ in the GNS
representation associated to $\tau $. Then define the $k$-th $L^{2}$-homology
group of $A$ to be the Hochschild homology group\[
H_{k}^{(2)}(A,\tau )=H_{k}(A;M\bar{\otimes }M^{o}).\]
Also define the $k$-th $L^{2}$-Betti number of $A$ to be its extended
Murray-von Neumann dimension (in the sense of Luck),\[
\beta _{k}^{(2)}(A,\tau )=\dim _{M\bar{\otimes }M^{o}}H_{k}(A;M\bar{\otimes }M^{o}).\]

\end{defn}
\n This definition of course depends on the trace
that we choose on $A$.

\subsection{The bar resolution.}

We will now give an {}``explicit'' description of the $L^{2}$-homology
group of $A$, which is at the same time the only description we have
at present. The situation is somewhat analogous to the case of dealing
with a general group $\Gamma $, for which one does not know whether
or not one can choose a {}``nice'' topological space on which $\Gamma $
can act properly and freely (which would give us a {}``nice'' resolution with
which to compute the homology). Thus one must instead resort to using
the universal cover $E\Gamma $ of the universal classifying space
$B\Gamma $ of $\Gamma $. 

\n Let $C_{k}(A)=(A\otimes A^{o})\otimes A^{\otimes k}$, $k=0,1,\ldots $,
viewed as a left $A\otimes A^{o}$-module via\[
(m \otimes n^{o})\cdot((a\otimes b^{o})\otimes a_{1}\otimes \cdots \otimes a_{k})=(m\,a\,\otimes (b\,n)^{o})\otimes
a_{1}\otimes \cdots \otimes a_{k},\]
where $m,n,a,b,a_{1},\ldots ,a_{k}\in A$. Note that
$C_{k}(A)$ is a free $A\otimes A^{o}$-module (typically on an infinite
number of generators if $k>0$). Define\[
\partial _{k}:C_{k}(A)\to C_{k-1}(A)\]
by the formula\begin{eqnarray}
\partial _{k}(T\otimes a_{1}\otimes \cdots \otimes a_{k}) & = & T\, (a_{1}\otimes 1)\, \otimes \cdots
\otimes a_{k}+\label{eq:differential}\\(-1)^{k}T\,(1 \otimes a_{k}^{o})\,\otimes a_{1}\otimes \cdots \otimes
a_{k-1}
 &+  & \sum _{j=1}^{k-1}(-1)^{j}T\otimes a_{1}\otimes \cdots \otimes
a_{j}a_{j+1}\otimes \cdots \otimes a_{k};\nonumber 
\end{eqnarray}
where $T\in A\otimes A^{o}$, $a_{1},\ldots ,a_{k}\in A$
and we use the algebra structure of $A\otimes A^{o}$
to define $T\,(a_{1}\otimes 1)$ and $T\,(1 \otimes a_{k}^{o})$. 
Then $(C_{*}(A),\partial _{*})$
is exact \cite{cartan-eilenberg:homalgebra}
and forms a resolution of the $A\otimes A^{o}$-left module
 $A$ of (\ref{tri}) by $A\otimes A^{o}$-left modules, with the last map $C_{0}(A)=A\otimes A^{o}\mapsto A $ given by the multiplication $m$, 
$$
\sum a_{i}\otimes b_{i}^{o} \rightarrow \sum a_{i}\, b_{i}
$$

\n Let $C_{k}^{(2)}(A)=M\bar{\otimes }M^{o}\otimes _{A\otimes A^{o}}C_{k}=(M\bar{\otimes
}M^{o})\otimes A^{\otimes k}$.
Let $\partial _{k}^{(2)}$ be given by the same formula as in (\ref{eq:differential}),
except that we now allow $T\in M\bar{\otimes }M^{o}$. Then
\begin{equation}
H_{k}^{(2)}(A,\tau )=\frac{\ker \partial _{k}^{(2)}}{\Image \partial _{k+1}^{(2)}}.
\end{equation}
is by construction a left $M\bar{\otimes
}M^{o}$-module.

\n It is in general not clear how to compute these homology groups (or their dimensions over
$M\bar{\otimes }M^{o}$). However, one can give a description of the
$L^{2}$-Betti numbers of $A$ as a  limit of finite numbers.

\begin{lem}
The bar resolution of $A$ can be written as an inductive limit of
its sub-complexes $C_{*}^{(\ell)}$, $\ell\in I$, so that each $C_{*}^{(\ell)}$
is a finitely generated left module over $A\otimes A^{o}$. Moreover, if
we denote by $i_{\ell,k}:C_{*}^{(\ell)}\to C_{*}^{(k)}$ the inclusion map,
and by $H_{*}{(\ell)}$ the homology of the complex $(M\bar{\otimes }M)\otimes
_{A\otimes A^{o}}C_{*}^{(\ell)}$,
then we have\begin{eqnarray}
H_{n}^{(2)}(A,\tau ) & = & {\lim_{\longrightarrow} }
(H_{n}{(\ell)},(i_{\ell,k})_{*}),\label{eq:homindlim}\\
\beta _{n}^{(2)}(A,\tau ) & = & \sup _{\ell}\inf _{k>\ell}\dim _{M\bar{\otimes
}M^{o}}(i_{\ell,k})_{*}H_{n}{(\ell)}.\label{eq:betaaslim}
\end{eqnarray}
 
\end{lem}
\begin{proof}
For each integer $n$, and  finite subset $F\subset A$, let $V=\Span F$,
and let $V_{0}=V^{\otimes n}$. Denote by $m$ the multiplication
map $m:A\otimes A\to A$, and let\begin{eqnarray*}
V_{1} & = & \Span ((m\otimes 1^{\otimes n-2})V_{0},(1\otimes m\otimes 1^{\otimes
n-3})V_{0},\ldots ,(1^{\otimes n-2}\otimes m)V_{0})\\
V_{2} & = & \Span ((m\otimes 1^{\otimes n-3})V_{1},(1\otimes m\otimes 1^{\otimes
n-4}V_{1},\ldots ,(1^{\otimes n-3}\otimes m)V_{1})
\end{eqnarray*}
and so on, where $1$ stands for the identity map $A \mapsto A$. Let\[
C_{k}^{(F,n)}=\begin{cases} (A\otimes A^{o})\otimes V_{n-k}, & k\leq n\\
 0, & k>n.\\ \end{cases}
\]
Note that $C_{k}^{(F,n)}$ forms a sub-complex of the bar resolution
of $A$.

\n Order the pairs $(F,n)$ by saying that $(F,n)<(F',n')$ if $F\subset F'$
and $n\leq n'$. If $(F,n)<(F',n')$, let\[
i_{(F,n)}:C_{k}^{(F,n)}\to C_{k}^{(F',n')}\]
be the inclusion map. Then each $C_{k}^{(F,n)}$ is finitely-generated
over $A\otimes A^{o}$ (it has at most $\dim V_{k}\,$ generators), and
the bar resolution of $A$ is the inductive limit of the sub-complexes
$C_{k}^{(F,n)}$. 

\n It follows that the homology group $H_{n}(A;M\bar{\otimes }M^{o})$
is itself the inductive limit of the directed system $\{H_{n}{(\ell)},i_{\ell,k}\}$,
so that (\ref{eq:homindlim}) holds.

\n Because of the finite generation assumptions, we have that\[
\dim _{M\bar{\otimes }M^{o}}(i_{\ell,k})_{*}H_{n}{(\ell)}<\infty ,\qquad \forall k>\ell.\]
Thus \cite[Theorem 2.9(2)]{luck:foundations1} implies (\ref{eq:betaaslim}).
\end{proof}
A reader interested in an even more explicit description, valid for
the first Betti number, is urged to look ahead to section \ref{sub:beta1inductive}.

\subsection{Group algebras.}

We have defined $L^{2}$-Betti numbers for any tracial $\ast$-algebra $(A,\tau )$.
While the case of interest is  when $A$ is a von Neumann algebra
(so that $M=A$), we want to point out that the definition works well
even in the purely algebraic setting.

\begin{prop}
Let $\Gamma $ be a discrete group, and denote by $\tau $ the von
Neumann trace on the group algebra $\mathbb{C}\Gamma $. Let $\Gamma^{(2)}=\Gamma \times \Gamma^{o} $, with $\Gamma^{o} $ the opposite group,
and view $\Gamma $ as a subgroup of $\Gamma^{(2)} $ via the diagonal inclusion
map $\g \mapsto \Delta(\g)=(\g, (\g^{-1})^{o})$. Then\begin{eqnarray}
H_{k}^{(2)}(\mathbb{C}\Gamma ,\tau ) & = & L(\Gamma^{(2)} )\otimes _{L(\Gamma
)}H_{k}^{(2)}(\Gamma ),\label{eq:algvsgroup1}\\
\beta _{k}^{(2)}(\mathbb{C}\Gamma ,\tau ) & = & \beta _{k}^{(2)}(\Gamma
).\label{eq:algvsgroup2}
\end{eqnarray}

\end{prop}
\begin{proof}
Note that $L(\Gamma )\bar{\otimes }L(\Gamma^{o} )=L(\Gamma \times \Gamma^{o} )=L(\Gamma^{(2)} )$.
By \cite{cartan-eilenberg:homalgebra}, we have that\[
H_{k}(\mathbb{C}\Gamma ;L(\Gamma )\bar{\otimes }L(\Gamma^{o} ))=H_{k}(\Gamma ,L(\Gamma^{(2)}
)).\]
Since $L(\Gamma^{(2)} )=L(\Gamma^{(2)} )\otimes _{L(\Gamma )}L(\Gamma )$ as $\Gamma $-modules,
and the functor $L(\Gamma^{(2)} )\otimes _{L(\Gamma )}$ is flat \cite[Theorem
3.3(1)]{luck:foundations1},
it follows that\[
H_{k}(\mathbb{C}\Gamma ;L(\Gamma )\bar{\otimes }L(\Gamma^{o} ))=L(\Gamma^{(2)} )\otimes
_{L(\Gamma )}H_{k}(\Gamma ;L(\Gamma ))=L(\Gamma^{(2)} )\otimes _{L(\Gamma
)}H_{k}^{(2)}(\Gamma ).\]
Equation (\ref{eq:algvsgroup2}) now follows from \cite[Theorem
3.3(2)]{luck:foundations1}.
\end{proof}

\subsection{Compressions of von Neumann algebras.}

As was shown by Gaboriau in the context of measurable measure-preserving
equivalence relations \cite{gaboriau:ell2}, $L^{2}$-Betti numbers
behave well under restrictions of equivalence relations. More precisely,
if an ergodic equivalence relation $R$ is restricted to subset $X$
of measure $\lambda $, then one has\[
\beta _{n}^{(2)}(R_{\lambda })=\frac{1}{\lambda }\beta _{n}^{(2)}(R).\]
The analogue of this fact is given by the following theorem. It should
be noted that the factor $1/\lambda $ in Gaboriau's result is replaced
in our context by the factor of $1/\lambda ^{2}$. This is explained
by the fact that in constructing our $L^{2}$-homology, we have passed
to the category of bimodules, so the natural object that we are working
with is $M\otimes M^{o}$ (and not $M$). Compressing $M$ to a projection
of trace $\lambda $ amounts to compressing $M\otimes M^{o}$ by a
projection of trace $\lambda ^{2}$.

\begin{thm}
Let $M$ be a factor and let $p\in M$ be a projection of trace $\lambda $.
Then\begin{equation}
\beta _{n}^{(2)}(pMp,\frac{1}{\tau (p)}\tau (p\cdot p))=\frac{1}{\lambda ^{2}}\beta
_{n}^{(2)}(M,\tau ).\label{eq:bettascale}\end{equation}

\end{thm}
\begin{proof}
Let $(C_{*}(M), \partial_{*})$ be the bar resolution of $M$,\begin{eqnarray*}
C_{k}(M) & = & (M\otimes M^{o})\otimes M^{\otimes k},\qquad k\geq 0
\end{eqnarray*}
with $\partial_{*}$ as in \ref{eq:differential}.
Let $N = M\otimes M^{o}, q=p\otimes p^{o}$ which is an idempotent in $N$. The 
reduced algebra $N_q= q\,N\,q$ is $pMp\otimes (pMp)^{o}$.

\n Let $F$ be the functor $V \mapsto F(V)=q\, V$
from the category of left $N$-modules to that of left $N_q$-modules.
It is an exact functor since for $T : V\mapsto W\, , \;x\in \Image T \cap q\, W$
one has $$x = T y = q T y = T q y \in \Image (T/qV)\,.$$  Note that in our case $qN$ is a projective left module over $N_q$ since $pM$ is a projective left module over $pMp$
and similarly for the opposite algebras.
Thus $F$ maps projective modules to projective modules.
Moreover 
when we apply $F$  to the ``trivial" $N$-module $M$ of \ref{tri}
we get $F(M)= pMp$, the ``trivial" $pMp\otimes (pMp)^{o}$-module.

\n This shows that $F(C_{*}(M), \partial_{*})$ is a projective resolution of 
$pMp$ and hence that,
\begin{eqnarray*}
H_{*}^{(2)}(pMp,\frac{1}{\tau (p)}\tau (p\cdot p)) & = & (p\otimes p^{o})H_{*}^{(2)}(M,\tau ).
\end{eqnarray*}

\n Equation (\ref{eq:bettascale}) now follows from the fact that for
an $M\bar{\otimes }M^{o}$ left module $V$ and the projection $q=p\otimes p^{o}\in
M\bar{\otimes }M^{o}$,
we have\[
\dim _{pMp\bar{\otimes (}pMp)^{o}}qV=\dim _{q(M\bar{\otimes
}M^{o})q}qV=\frac{1}{(\tau \otimes \tau )(q)}\dim _{M\bar{\otimes }M^{o}}V.\]

\end{proof}

\subsection{Direct sums.}

$L^{2}$-homology behaves well with respect to direct sums:

\begin{prop}
Let $(A,\tau )=\bigoplus _{i}(A_{i},\tau _{i})$ (finite direct sum),
so that the trace on $A$ decomposes as $\bigoplus _{i}\alpha _{i}\tau _{i}$
in terms of normalized traces $\tau _{i}$ on $A_{i}$.

Then\[
H_{k}^{(2)}(A,\tau )\cong \bigoplus _{i}H_{k}^{(2)}(A_{i},\tau _{i})\]
 and\[
\beta _{k}^{(2)}(A,\tau )=\sum _{i}\alpha _{i}^{2}\beta _{k}^{(2)}(A_{i},\tau _{i}).\]

\end{prop}
\begin{proof}
Let $C_{\ast}^{(i)}$
be the bar resolution of $A_{i}$ with its differential $d_{\ast}^{(i)}$,
and let\[
C_{k}=\bigoplus _{i}C_{k}^{(i)},\qquad d_{k}=\bigoplus _{i}d_{k}^{(i)}.\]
 Then each $C_{k}$, $k\geq 1$ is a projective module over $A\otimes A^{o}$.
This is because $\bigoplus _{i}A_{i}\otimes A_{i}^{o}$ is a direct
summand of $A\otimes A^{o}=\bigoplus _{i,j}A_{i}\otimes A_{j}^{o}$.
Thus $(C_{k},d_{k})$ is a projective resolution of $A$. Using
this resolution to compute the $L^{2}$-homology of $A$ we obtain
\[
H_{k}^{(2)}(A,\tau )=\bigoplus _{i}H_{k}^{(2)}(A_{i},\tau _{i}).\]

\n Let $M_{i}=W^{*}(A_{i})$, $M=W^{*}(A)$ (each time in the GNS representation
associated to $\tau $). The formula for Betti numbers is now immediate,
once we remark that if $V_{i}$ is a module over $M_{i}\bar{\otimes }M_{i}^{o}$,
then\[
\dim _{M\bar{\otimes }M^{o}}(\bigoplus _{i}V_{i})=\sum \alpha _{i}^{2}\dim
_{M_{i}\bar{\otimes }M_{i}^{o}}V_{i},\]
 the factor $\alpha _{i}^{2}$ coming from the fact that\[
M\bar{\otimes }M^{o}=\bigoplus _{i,j}M_{i}\otimes M_{j}^{o},\qquad \tau \otimes \tau
=\bigoplus _{i,j}\alpha _{i}\alpha _{j}\tau _{i}\otimes \tau _{j}.\]

\end{proof}

\subsection{Zeroth Betti number and zeroth $L^{2}$-homology for  von Neumann
algebras.}

Let $M$ be a von Neumann algebra. By definition $$H_{0}^{(2)}(M,\tau )
=\frac{\ker \partial _{0}^{(2)}}{\Image \partial _{1}^{(2)}}$$
is  the quotient
of $M\bar{\otimes }M^{o}$ by the left ideal $L$ generated by \[
V=\{1\otimes n^{o}-n\otimes 1:n\in M\}.\]
or in other words \[
H_{0}^{(2)}(M,\tau )=(M\bar{\otimes }M^{o})\otimes _{M\otimes M^{o}}M.\]

\begin{prop}
Let $M$ be a II$_1$-factor. Then $H_{0}^{(2)}(M,\tau )\neq 0$
if and only if $M$ is hyperfinite. 
\end{prop}
\begin{proof}
 Since as a left $M\bar{\otimes }M^{o}$-module, 
 $H_{0}^{(2)}(M)$ is clearly generated by the class in the
quotient of the element $1\otimes 1$, $H_{0}^{(2)}(M)=0$ if and
only if $[1\otimes 1]=0$.

\n By (\cite{connes:injective})  $M$ is
hyperfinite if and only if the trivial correspondence is weakly contained
in the coarse correspondence. That is to say, there exists a nonzero
non-normal state $\theta :M\bar{\otimes }M^{o}\rightarrow \mathbb{C}$
with the property that
\begin{equation} \label{hyper}
\theta (m\otimes n^{o})=\tau (m\, n),\qquad \forall m,n\in M.
\end{equation}

\n Assume first that $M$ is hyperfinite. Then
$\theta(\,x^\ast \,x)=0$ for $x\in V$. Thus $\theta |_{J}=0$, so
that $\theta $ factors through to a non-zero linear functional on
$H_{0}^{(2)}(M)$ and $H_{0}^{(2)}(M)\neq 0$.

\n Conversely, assume that $H_{0}^{(2)}(M)\neq 0$. Then for any $n$ unitaries
$u_{i}\in M$, the operator \[
T=\frac{1}{n}\sum u_{i}\otimes u_{i}^{*o}-1,\]

 belongs to the left ideal $L$ so that $T$ and hence $T^{*}T$ are
not invertible, for any $u_{i}$. Thus one can find a non-normal state
$\phi $ on $M\bar{\otimes }M^{o}$, for which $\phi (T^{*}T)=0$ for
any such $T$. Denoting by $\xi _{\phi }$ the associated cyclic vector,
we get that $T\xi _{\phi }=0$ and so\[
(u_{i}\otimes u_{i}^{*o})\,\xi _{\phi }=\xi _{\phi }.\]
 But then we have\[
(u_{i}\otimes 1)\,\xi _{\phi }=(1\otimes u_{i}^{o})\,\xi _{\phi }\]
 and hence\[
(n\otimes 1)\,\xi _{\phi }=(1\otimes n^{o})\,\xi _{\phi }\]
 for all $n\in M$. Hence\[
(m\,n\otimes 1)\,\xi _{\phi }=(m \otimes 1)(n\otimes 1)\,\xi _{\phi }=(m\otimes n^{o})\,\xi _{\phi },\]
 so that $\phi (m\otimes n^{o})=\phi (m\,n\otimes 1)$. Lastly,\begin{align*}
\phi (n\,m\otimes 1) &=\langle (n\otimes m^{o})\,\xi _{\phi },\xi _{\phi }\rangle \\
 & =\langle (1\otimes m^{o})\cdot (n\otimes 1)\,\xi _{\phi },\xi _{\phi }\rangle \\
 & =\langle (n\otimes 1)\,\xi _{\phi },(1\otimes m^{*o})\xi _{\phi }\rangle \\
 & =\langle (n\otimes 1)\,\xi _{\phi },(m^{*}\otimes 1)\,\xi _{\phi }\rangle \\
 & =\langle (m\, n\otimes 1)\,\xi _{\phi },\xi _{\phi }\rangle = \phi (m\,n\otimes 1),
\end{align*} so that $\phi |_{M\otimes 1}$ is a trace, and hence the unique trace
$\tau $ on $M\cong M\otimes 1$, thus $\phi$ fulfills \ref{hyper}. 
\end{proof}
We get as consequence of Luck's Theorem 1.8 \cite{luck:foundations1},
his definition of the projective part of a module and of Theorem 0.6
in \cite{luck:foundations1}, and of the fact that $H_{0}^{(2)}(A,\tau
)=M\bar{\otimes }M^{o}\otimes _{A\otimes A^{o}}A$,
$M=W^{*}(A)$ is finitely (in fact, singly) generated as an $M\bar{\otimes }M^{o}$
module, that\[
\beta _{0}(A,\tau )=\dim _{M\bar{\otimes }M^{o}}\textrm{Hom}(M\bar{\otimes
}M^{o}\otimes _{A\otimes A^{o}}A,M\bar{\otimes }M^{o})\]

\begin{prop}
\label{pro:centerbetti0}If $(A,\tau )$ contains an element $x$,
whose distribution with respect to $\tau $ is non-atomic, then $\beta
_{0}^{(2)}(A,\tau )=0$. 
\end{prop}
\begin{proof}
Assume $\beta _{0}^{(2)}(A,\tau )\neq 0$, let $\phi \neq 0$, 
$\phi \in \Hom (M\bar{\otimes }M^{o}\otimes _{A\otimes A^{o}}A,M\bar{\otimes }M^{o})$.
Denote by $[1\otimes 1]$ the class of $1\otimes 1$ in $M\bar{\otimes }M^{o}\otimes
_{A\otimes A^{o}}A$.
Let $\xi =\phi ([1\otimes 1])\in M\bar{\otimes }M^{o}$. Since
$M\bar{\otimes }M^{o}\otimes _{A\otimes A^{o}}A$ is generated by
$[1\otimes 1]$, $\phi \neq 0$ implies that $\xi \neq 0$.

\n We thus have a vector $\xi \neq 0$ in $L^{2}(M\bar{\otimes }M^{o})$
with the property that $(m\otimes 1-1\otimes m^{o})\,\xi =0$ for all $m\in A$ and hence
for all $m\in A''=M$. Identify $L^{2}(M\bar{\otimes }M^{o})$ with the
space of Hilbert-Schmidt operators on $L^{2}(M)$, by the 
map $\Psi$ of (\ref{bim}).
 Then $\Psi(\xi) $ is a non-zero Hilbert-Schmidt operator, commuting
with $M$ by (\ref{bim}). But this is impossible, since $M$ contains an element
with a diffuse spectrum.
\end{proof}
\begin{cor}
If $M$ is a type II$_{1}$ factor, then $\beta _{0}^{(2)}(M)=0$. 
\end{cor}
\begin{prop}
Let $M$ be a finite-dimensional von Neumann algebra with a faithful
trace $\tau $. Decompose $M=\oplus M_{i}$ into factors with $M_{i}\cong
M_{k_{i}\times k_{i}}$(the
algebra of $k_{i}\times k_{i}$ matrices), and let $\lambda _{i}$
be the trace of the minimal central projection in $M$ corresponding
to the $i$-th summand.

Then \begin{eqnarray*}
\beta _{0}^{(2)}(M) & = & \sum _{i}\frac{\lambda _{i}^{2}}{k_{i}^{2}},\\
\beta _{k}^{(2)}(M) & = & 0,\qquad k\geq 1.
\end{eqnarray*}

\end{prop}
\begin{proof}
Since there is no difference between $M\otimes M$ and $M\bar{\otimes }M$
in the finite-dimensional case, $\beta _{k}^{(2)}(M)=0$ if $k>0$.
For $\beta _{0}^{(2)}$ we find easily that $\beta _{0}^{(2)}(\mathbb{C})=1$;
the compression formula then gives $\beta _{0}^{(2)}(M_{k\times k})=\frac{1}{k^{2}}$,
and the direct sum formula gives us finally the desired expression. 
\end{proof}

\subsection{$L^{2}$-Betti numbers for bimodule maps.}

\subsubsection{Betti numbers for group module maps.}

$L^{2}$-Betti numbers can be more generally defined for maps between
group modules. Let us for definiteness consider a module map $f$
between two free left $\Gamma$-modules:\[
f:(\mathbb{C}\Gamma )^{n}\to (\mathbb{C}\Gamma )^{m}.\]
Thus $f$ is given by an $n\times m$ matrix of right-multiplication
operators in $\mathbb{C}\Gamma $. Consider now\[
f^{(p)}:(\ell ^{p}(\Gamma ))^n\to (\ell ^{p}(\Gamma ))^m,\]
given by the same matrix. The kernel of $f^{(p)}$ may be larger than
the $\ell ^{p}$-closure of the kernel of $f$. To measure the difference,
consider for $p\leq 2$\[
\beta ^{(2,p)}(f)=\dim _{L(\Gamma )}\frac{\overline{\ker f^{(p)}}^{\ell
^{2}}}{\overline{\ker f}^{\ell ^{2}}}.\]
Set\[
\beta ^{(2)}(f)=\beta ^{(2,2)}(f)\]
for convenience. 

\n Note that if $\Gamma $ acts freely and cocompactly on some contractible 
simplicial complex $X$,
then $\beta _{j}^{(2)}(\Gamma )$ is then just $\beta ^{(2)}(\partial _{j})$,
where $\partial _{j}$ is the boundary operator of $X$. Indeed, contractibility implies that $\ker \partial _{k}=\Image \partial
_{k+1}$,
so that the closures of these two subspaces of $(\ell ^{2})^N$ are the same.

\subsubsection{Betti numbers for bimodule maps.}

Let $(A,\tau )$ be a tracial $\ast$-algebra. Let $F:(A\otimes A^{o})^{n}\rightarrow
(A\otimes A^{o})^{m}$
be a left $A\otimes A^{o}$-module map (or, equivalently, an $A,A$-bimodule
map). Then $F$ is given by a matrix\[
F=\left(\begin{array}{ccc}
 F_{11} & \cdots  & F_{1n}\\
 \vdots  & \ddots  & \vdots \\
 F_{m1} & \cdots  & F_{mn}\\ \end{array}
\right),\]
 where  $F_{ij}\in A\otimes A^{o}$ and the action
 is given by right multiplication. 

\n Let $M=W^{*}(A)$ in the GNS representation associated to $\tau $.

\n The right multiplication by $F_{ij}$ admits a unique continuous extension
to $L^{2}(M\bar{\otimes }M^{o})$. Thus $F$ admits a unique continuous extension to a left $M\bar{\otimes }M^{o}$-module map
from $(L^{2}(M\bar{\otimes }M^{o}))^{n}$ to $(L^{2}(M\bar{\otimes }M^{o}))^{m}$,
which we denote by $F^{(2)}$.

\n We note that $\ker F\subset \overline{\ker F}\subset \ker F^{(2)}$
(where $\overline{\cdot }$ refers to closure in $L^{2}$-norm). By
analogy with the group case, we make the following definition (compare
\cite{luck:hilbertmodules}, Definition 5.1).

\begin{defn}
The $L^{2}$ Betti number of $F$ is the Murray-von Neumann dimension\[
\beta ^{(2)}(F)=\dim _{M\bar{\otimes }M^{o}}\frac{\ker F^{(2)}}{\overline{\ker
F}}=\dim _{M\bar{\otimes }M^{o}}\ker F^{(2)}-\dim _{M\bar{\otimes
}M^{o}}\overline{\ker F}.\]

\end{defn}
\begin{lem}
\label{lem:l2closure} One has \[
\dim _{M\bar{\otimes }M^{o}}\frac{\ker F^{(2)}}{\overline{\ker F}^{L^{2}}}=\dim
_{M\bar{\otimes }M^{o}}\frac{\ker F^{vN}}{(M\bar{\otimes }M^{o})\ker F},\]
 where $F^{vN}$ is the restriction of $F^{(2)}$ to $(M\bar{\otimes
}M^{o})^{n}\subset (L^{2}(M\bar{\otimes }M^{o}))^{n}$,
and $(M\bar{\otimes }M^{o})\cdot \ker F$ denotes the saturation of
$\ker F$ under the action of $M\bar{\otimes }M^{o}$. 
\end{lem}
\n The proof is almost verbatim the argument on the bottom of page 158
and top of page 159 of \cite{luck:foundations1}, see also Theorem
5.4 of \cite{luck:hilbertmodules}. 

\begin{rem}
Note that $F^{vN}$ is exactly the map\[
1\otimes F:(M\bar{\otimes }M^{o})\otimes _{A\otimes A^{o}}(A\otimes
A^{o})^{n}\rightarrow (M\bar{\otimes }M^{o})\otimes _{A\otimes A^{o}}(A\otimes
A^{o})^{m},\]
if we identify $(M\bar{\otimes }M^{o})\otimes _{A\otimes A^{o}}(A\otimes A^{o})$
with $M\bar{\otimes }M^{o}$. Thus in particular,\[
\beta ^{(2)}(F)=\dim _{M\bar{\otimes }M^{o}}\frac{\ker (1\otimes F)}{(M\bar{\otimes
}M^{o})\cdot \ker F}.\]

\end{rem}

\subsubsection{Homological algebra interpretation.}

Let $F:(A\otimes A^{o})^{n}\to (A\otimes A^{o})^{m}$ be a bimodule
map, as above. Put $V=(A\otimes A^{o})^{n}$, $W=(A\otimes A^{o})^{m}$.

\n Consider the exact sequence\[
V\stackrel{F}{\longrightarrow }W\to \Image F\to 0.\]
Since the $A\otimes A^{o}$- left 
modules $V$ and $W$ are projective (in fact, free), this
sequence is the beginning of a projective resolution of $\Image F$.
More precisely, one can choose projective modules $V_{1},V_{2},\ldots $
and morphisms $F_{1},F_{2},\ldots$ so that the following sequence is exact:\[
\cdots \to V_{2}\stackrel{F_{2}}{\longrightarrow
}V_{1}\stackrel{F_{1}}{\longrightarrow }V\stackrel{F}{\longrightarrow }W\to \Image
F\to 0.\]
Note that $\Image F_{1}=\ker F$. Consider the differential complex\[
\cdots \to (M\bar{\otimes }M^{o})\otimes _{A\otimes A^{o}}V_{1}\stackrel{1\otimes
F_{1}}{\longrightarrow }(M\bar{\otimes }M^{o})\otimes _{A\otimes
A^{o}}V\stackrel{1\otimes F}{\longrightarrow }(M\bar{\otimes }M^{o})\otimes
_{A\otimes A^{o}}W\to \cdots .\]
By definition,\[
\Tor _{1}^{A\otimes A^{o}}(\Image F,M\bar{\otimes }M^{o})=\frac{\ker (1\otimes
F)}{\Image (1\otimes F_{1})}.\]
Since $\Image (1\otimes F_{1})=(M\bar{\otimes }M^{o})\cdot \Image
F_{1}=(M\bar{\otimes }M^{o})\ker F$,
we conclude that
\begin{equation}\label{hom1}
\beta ^{(2)}(F)=\dim _{M\bar{\otimes }M^{o}}(\Tor _{1}^{A\otimes A^{o}}(\Image
F,M\bar{\otimes }M^{o})).
\end{equation}

\subsubsection{Examples of Betti numbers.}

The following statement gives many examples of bimodule maps over
von Neumann algebras, for which the $L^{2}$ Betti numbers are non-zero.
For a group module map $f:\mathbb{C}\Gamma ^{n}\rightarrow \mathbb{C}\Gamma ^{m}$,
denote by $f^{(1)}$ its extension to $\ell ^{1}(\Gamma )^{n}$. Denote
by $\beta ^{(2,1)}(f)$ the dimension \[
\beta ^{(2,1)}(f)=\dim _{L(\Gamma )}(\ker f^{(2)})-\dim _{L(\Gamma )}\overline{\ker
f^{(1)}}.\]
\medskip

\begin{thm}
\label{Prop:sameAsInGroupCase}Let $\Gamma $ be a discrete group,
$n,m<\infty $ and let $f:\mathbb{C}\Gamma ^{n}\rightarrow \mathbb{C}\Gamma ^{m}$
be a $\Gamma$-left module map given by a matrix with entries $f_{ij}\in \mathbb{C}\Gamma $.
Let $A=M=L(\Gamma )$. 

\n Let $F_{ij}=\Delta(f_{ij})\in M\bar{\otimes }M^{o}$ be the images of $f_{ij}$
under the canonical diagonal inclusion $\Delta(\g)=(\g, (\g^{-1})^{o})$ of $\mathbb{C}\Gamma $
into $M\otimes M^{o}=L(\Gamma )\otimes L(\Gamma^{o} )$. Let $F:(M\otimes
M^{o})^{n}\rightarrow (M\otimes M^{o})^{m}$
be given by the matrix whose entries are right multiplications by
$F_{ij}$. Then,

$$\qquad \beta ^{(2)}(f)\geq \beta ^{(2)}(F)\geq \beta ^{(2,1)}(f)$$ 
\end{thm}

\bigskip
\n Note that the statement is not automatic, since we are comparing\[
\Tor _{1}^{\Gamma }(\Image f,L(\Gamma ))\qquad \textrm{with}\qquad \Tor
_{1}^{L(\Gamma )\otimes L(\Gamma^{o} )}(\Image F,L(\Gamma )\bar{\otimes }L(\Gamma^{o} ))\]
and not\[
\Tor _{1}^{\Gamma }(\Image f,L(\Gamma ))\qquad \textrm{with}\qquad \Tor
_{1}^{\mathbb{C}\Gamma \otimes \mathbb{C}\Gamma^{o} }(\Image F,L(\Gamma )\bar{\otimes
}L(\Gamma^{o} ))\]

\n Let $\Delta_*$ be the induction functor from left $\Gamma$-modules
to left $M\otimes M^{o}$-modules associated to the morphism $\Delta:\mathbb{C}\Gamma 
\mapsto M\otimes M^{o}$,
$$
\Delta_*(X):= (M\otimes M^{o})\otimes _{\mathbb{C}\Gamma }X
$$
where $M\otimes M^{o}$ is viewed as a right $\mathbb{C}\Gamma $-module
using $\Delta$, then the $M\otimes M^{o}$-module
$(M\otimes M^{o})^{n}$ is induced from $\mathbb{C}\Gamma ^{n}$ while
$F$ is induced from $f$, or in short $F=1\otimes f.$\\

\begin{proof}
It is sufficient to prove that $\dim _{M\bar{\otimes }M^{o}}\ker F^{(2)}=\dim
_{L(\Gamma )}\ker f^{(2)}$,
$\dim _{M\bar{\otimes }M^{o}}\overline{\ker F}\geq \dim _{L(\Gamma )}\overline{\ker f}$
and $\dim _{M\bar{\otimes }M^{o}}\overline{\ker F}\leq \dim _{L(\Gamma )}\overline{\ker
f^{(1)}}$.
Here $f^{(1)}$ stands for the extension of $f$ to $\ell ^{1}(\Gamma )$.

\n The morphism $\Delta$ preserves the trace
and extends to an inclusion of von Neumann algebras $\Delta:
L(\Gamma )\mapsto  M\bar{\otimes }M^{o}$. Denote by $F^{*}$ the {}``adjoint'' of $F$ (with $i,j$-th entry
$F_{ji}^{*}$). Then $\ker F^{(2)}=\ker ((F^{(2)})^{*}F^{(2)})$.
Regard $T=(F^{(2)})^{*}F^{(2)}$ as an element of the algebra of $n\times n$
matrices over $M\bar{\otimes }M^{o}$. Then the dimension of the kernel
of $F^{(2)}$ is precisely the non-normalized trace, computed in $M_{n\times
n}(M\bar{\otimes }M^{o})$,
of the spectral projection $P$ of $T$ corresponding to the eigenvalue
$0$. But $F= \Delta(f)$, and hence $P= \Delta(p)$ where $p$
is the spectral projection $p$ corresponding to the eigenvalue
$0$ of the element $t=(f^{(2)})^{*}f^{(2)}\in M_{n\times n}(L(\Gamma ))$.
The trace of $p$ is exactly $\dim _{L(\Gamma )}\ker f^{(2)}$. Thus $$\dim _{M\bar{\otimes
}M^{o}}\ker F^{(2)}=\dim _{L(\Gamma )}\ker f^{(2)}$$

\n Consider now the orthogonal projection $E:L^{2}(M\bar{\otimes }M^{o})\rightarrow
\Delta(L^{2}(M))$. It defines a conditional expectation $E: M\bar{\otimes }M^{o}
\rightarrow M$ where we identify $\Delta M$ with $M$.
Note that if $\eta\in M\otimes M^{o}\subset L^{2}(M\bar{\otimes }M^{o})$,
then $E(\eta)\in \ell ^{1}(\Gamma )$. To prove this, note first that
it is sufficient to prove that $E(\eta)\in \ell ^{1}(\Gamma )$ whenever
$\eta$ is a simple tensor of the form $\xi \otimes \zeta $, $\xi \in L^{2}(M)$,
$\zeta \in L^{2}(M^{o})$. For $\gamma \in \Gamma $, let $u_{\gamma }$
be the corresponding unitary in $M$. Let $\xi =\sum \alpha _{\gamma }u_{\gamma }$
and $\zeta =\sum \beta _{\gamma }u_{\gamma ^{-1}}^{o}$ (where the sums
are in $L^{2}$ sense). Then $E(\xi \otimes \zeta )=\sum \alpha _{\gamma }\,\beta
_{\gamma }\,\Delta(u_{\gamma })\sim \sum \alpha _{\gamma }\,\beta
_{\gamma }\,u_{\gamma }$.
Since the sequences $\{\alpha _{\gamma }\}$ and $\{\beta _{\gamma }\}$
are in $\ell ^{2}(\Gamma )$, their product lies in $\ell ^{1}(\Gamma )$,
and $E(\eta)\in \ell ^{1}(\Gamma )$.

\n Denote by $E_{\gamma }$ the map $\eta\mapsto E((u_{\gamma ^{-1} }\otimes 1)\cdot \eta$).
Then any $\eta\in L^{2}(M\bar{\otimes }M^{o})$ is the $L^{2}$-sum $\eta=\sum
(u_{\gamma }\otimes 1)\cdot E_{\gamma
}(\eta)$.
We extend $E$ and $E_{\gamma }$ (componentwise) to maps of direct
sums of $L^{2}(M\bar{\otimes }M^{o})$ and denote them with the same
letter. Since $E$ is a conditional expectation and $F$ acts by right multiplication
by the $\Delta(f_{ij})$ one has $E\circ F^{(2)}=f^{(2)}\circ E$ and also $$E_{\gamma }\circ
F^{(2)}=f^{(2)}\circ E_{\gamma }\, \quad \forall \gamma \in \Gamma $$

\n Let now $\eta\in \ker F$. Then $\eta\in (M\otimes M^{o})^n$, and $E_{\gamma }(\eta)\in (\ell
^{1}(\Gamma ))^{n}$,
for all $\gamma $. Denote by $f^{(1)}$ the restriction of $f^{(2)}$
to $(\ell ^{1}(\Gamma ))^{n}$. Then $f^{(1)}\circ E_{\gamma }(\eta)
=E_{\gamma }\circ F^{(2)}(\eta)=0$,
so that $E_{\gamma }(\eta)\in \ker f^{(1)}$. Denote by $q$ the projection
in $M_{n\times n}(M)$ corresponding to the invariant subspace $\overline{\ker
f^{(1)}}$.
Let $Q =\Delta(q)$  be the image of $q$ 
under the inclusion map of $M_{n\times n}(M)\subset
M_{n\times n}(M\bar{\otimes }M^{o})$
induced by $\Delta$.
Since $E_{\gamma }(\eta)\in \ker f^{(1)}$ for all $\gamma $, it follows
easily that $\eta=\sum
(u_{\gamma }\otimes 1)\cdot E_{\gamma
}(\eta)$
is in the range of $Q$. Thus we have proved that $\ker F$ is contained
in the range of $Q$, so that $\dim _{M\bar{\otimes }M^{o}}\overline{\ker F}\leq \dim
_{M\bar{\otimes }M^{o}}\Image Q=\dim _{M}\Image q=\dim _{M}\overline{\ker f^{(1)}}$.\\

\n Finally, if we are given a finite sequence $t_{\gamma }\in \ker f$,
then $\sum
(u_{\gamma }\otimes 1)\cdot t_{\gamma }\in \ker F$.
Thus the induced module from $\overline{\ker f}$
is contained in $\overline{\ker F}$.
This shows that $\dim _{M}\overline{\ker f}\leq \dim
_{M\bar{\otimes }M^{o}}\overline{\ker F}$.
Thus we have\[
\dim _{M}\overline{\ker f}\leq \dim _{M\bar{\otimes }M^{o}}\overline{\ker F}\leq
\dim _{M}\overline{\ker f^{(1)}}.\]

\end{proof}
It would be interesting to know exactly when $\beta ^{(2)}(F)=\beta ^{(2)}(f)$.
Note that by the results of \cite{mineev:vanishing}, if $\Gamma $
is a {}``combable group'' (in particular, a finitely generated hyperbolic
group), and $f:C_{n}(X)\rightarrow C_{n-1}(X)$ is the boundary homomorphism
of a contractible chain complex with a cocompact free action
of $\Gamma $, then $\beta ^{(2)}(f)=\beta ^{(2,1)}(f)$. Indeed,
it it proved in \cite{mineev:vanishing} that any element in $\ker f^{(1)}$
can be approximated in $\ell ^{1}$ (and hence $\ell ^{2}$) by elements
from $\ker f$. This implies the following fact:

\begin{thm}
\label{thm:combableBetti}Let $\Gamma $ be a discrete combable group
acting freely and co-compactly on a contractible chain complex $C_{*}(X)$
of a topological space $X$. Let $D_{*}=(L(\Gamma )\otimes L(\Gamma^{o} ))\otimes
_{\Gamma }C_{*}(X)$,
where $\Gamma $ acts on $L(\Gamma )\otimes L(\Gamma^{o} )$ by $\Delta$.
Then for each $k$,\begin{eqnarray*}
\dim _{L(\Gamma )\bar{\otimes }L(\Gamma^{o} )}H_{k}((L(\Gamma )\bar{\otimes }L(\Gamma^{o}
))\otimes _{L(\Gamma )\otimes L(\Gamma^{o} )}D_{*}(X)) & = & \dim _{L(\Gamma
)}H_{k}(L(\Gamma )\otimes _{\Gamma }C_{*}(X))\\
 & = & \beta _{k}^{(2)}(\Gamma ).
\end{eqnarray*}

\end{thm}
\n This theorem is of interest in conjunction with equation (\ref{eq:betaaslim}),
since $D_{*}(X)$ (is homotopic to a a sub-complex that) occurs among
the approximating sub-complexes of the bar resolution of $L(\Gamma )$.

\subsection{Dual definition of Betti numbers for bimodule maps.}

For the remainder of this section, we shall concentrate on bimodule
maps over a von Neumann algebra; i.e., we assume that $A=M$ is a
von Neumann algebra with a fixed trace $\tau $.

\n Let $f:(M\otimes M^{o})^{n}\rightarrow (M\otimes M^{o})^{m}$ be a
left $M\otimes M^{o}$-module map. We are interested in the size of the kernel $\ker f$
in $(M\otimes M^{o})^{n}$. As before, denote by $f^{vN}$ and $f^{(2)}$,
respectively, the extensions of $f$ to $(M\bar{\otimes }M^{o})^{n}$ and
$L^{2}(M\bar{\otimes }M^{o})^{n}$.

 \n Let us consider the algebraic tensor product $FR=L^{2}(M)\otimes L^{2}(M^{o})$
as a subset of $L^{2}(M)\bar{\otimes }L^{2}(M^{o})=L^{2}(M\bar{\otimes }M^{o})$
in the natural way. Note that $M\otimes M^{o}\subset FR$.

\n Note also that in the identification $\Psi $ of $L^{2}(M\bar{\otimes }M^{o})$
with the space $HS=HS(L^{2}(M))$ of Hilbert-Schmidt operators on
$L^{2}(M)$ (see (\ref{eq:defofPsi})), the set $FR$ corresponds
precisely to the subset of $HS$ consisting of finite-rank operators.

\n We begin with a lemma, which shows that it does not matter for the
purposes of $L^{2}$-closure whether we compute the kernel of $f^{(2)}$
in $(M\otimes M^{o})^{n}$ or $FR^{n}$.

\begin{lem}
\label{lem:algsameasFR}$\ker f^{(2)}\cap FR^{n}$ has the same $L^{2}$-closure
in $L^{2}(M\bar{\otimes }M^{o})^{n}$ as $\ker f=\ker f^{(2)}\cap (M\otimes M^{o})^{n}$.
Thus\begin{equation}
\beta ^{(2)}(f)=\dim _{M\bar{\otimes }M^{o}}\ker f^{(2)}-\dim _{M\bar{\otimes }M^{o}}(\overline{\ker f^{(2)}\cap FR^{n}})\label{eq:betaFR}\end{equation}

\end{lem}
\begin{proof}
We view elements of $L^{2}(M)$ as
 unbounded operators (of left multiplication) 
on $L^{2}(M)$, affiliated with $M$.

\n Given any finite subset $K \subset L^{2}(M)$, and $\epsilon >0$
there exists a projection $e \in M$ such that 
$$
e \,\xi \in M \,,\quad \Vert e \xi - \xi \Vert < \epsilon \,,\quad \forall \xi \in K 
$$
where $e \xi \in M$ means that the a priori unbounded operator 
of left multiplication by $e \xi $
is bounded. This is proved for a single $\xi$ using the polar decomposition $\xi=bu$
and a suitable spectral projection of the unbounded self-adjoint operator $b$.
One controls moreover the trace $\tau(1-e)<\epsilon$. Taking the intersection
$f$ of the 
projections $e_\xi, \xi \in K$ one gets $\tau(1-f)<n \epsilon$,
$n= {\rm card}(K)$. This gives a sequence of projections such that 
$f_k \,\xi$ is bounded $\forall \xi \in K $ and $f_k \rightarrow 1$
in $L^2$ and hence strongly, which gives the answer. This shows that for
any element $\xi \in FR^{n}$ there exists a projection $e \in M$
such that $(e \otimes e^{o})\,\xi \in (M\otimes M^{o})^{n}$
and $ \Vert (e \otimes e^{o})\,\xi - \xi \Vert < \epsilon$. Since 
$\ker f^{(2)}\cap FR^{n}$ is  a left module over $M\otimes M^{o}$
the conclusion of the lemma follows.
\end{proof}
As consequence, we have the following description of the dimension
of the kernel of a left $M\otimes M^{o}$-module
 map $f:(M\otimes M^{o})^{n}\rightarrow (M\otimes M^{o})^{m}$.
Extend $f$ to a map (still denoted by $f$) from $FR^{n}\rightarrow FR^{m}$
as in the previous Lemma. Identify $FR$ with the set of finite-rank
operators on $L^{2}(M)$ using the identification $\Psi $. 
Let $B$ be the von Neumann algebra of all
bounded operators on $L^{2}(M)$. Denote by $\langle \cdot ,\cdot \rangle $
the canonical pairing between $FR^{n}$ and $B^{n}$ given by \[
\langle (T_{1},\ldots ,T_{n}),(S_{1},\ldots ,S_{n})\rangle =\sum _{j=1}^{n}\Tr (T_{j}S_{j}).\]
 Denote by $f^{\,t}$ the map \[
f^{\,t}:B^{m}\rightarrow B^{n}\]
 uniquely determined by \[
\langle f^{\,t}(T),S\rangle =\langle T,f(S)\rangle .\]
 The map $f^{\,t}$ is well-defined and is in fact given by 
right multiplication by the matrix $f_{ji}^{\,t}$
where for $h =\sum m_i \otimes n_i^{o} \in M\otimes M^{o}$ one lets
$$
h^{\,t}:= \sum n_i \otimes m_i^{o} \,.
$$
The right multiplication by $a\otimes b^{o}\in M\otimes M^{o}$ acts
in the obvious way
 on $(M\otimes M^{o})^{n}$ and becomes, after the identification $\Psi $, using (\ref{bim})
\begin{equation}\label{rm}
T \cdot (a\otimes b^{o})=Ja^{*}J\,T\,Jb^{*}J \,,\quad \forall T \in B\,.
\end{equation}

\medskip

\begin{lem}
\label{lem:dualbetti} We have \[
\dim _{M\bar{\otimes }M^{o}}((M\bar{\otimes }M^{o})\cdot \ker f)=n-\dim _{M\bar{\otimes }M^{o}}(\overline{f^{\,t}(B^{m})}^{w}\cap HS^{n}),\]
 where the closure is taken with respect to the weak operator topology.
In particular, \[
\beta ^{(2)}(f)=\dim _{M\bar{\otimes }M^{o}}\left(\frac{\overline{f^{\,t}(HS^{m})}^{w}\cap HS^{n}}{\overline{f^{\,t}(HS^{m})}^{HS}}\right).\]

\end{lem}
\begin{proof}
We note that $FR$ is the dual of $B(H)$ with respect to the weak
topology. By duality, $\overline{f^{\,t}(B^{m})}^{w}$ is the annihilator of $\ker f\subset FR^{n}$. Let us show that $\overline{f^{\,t}(B^{m})}^{w}\cap HS$
is the annihilator of $\overline{\ker f}^{HS}$ in $HS^{n}$.
The answer then follows by duality in  $HS^{n}$ whose dimension over 
$M\bar{\otimes }M^{o}$ is equal to $n$. 

\n Note that the two pairings are compatible. By continuity of the pairing
in $HS^{n}$, 
$\overline{f^{\,t}(B^{m})}^{w}\cap HS^{n}$ is perpendicular
to $\overline{\ker f}^{HS}$.

\n Since the Hilbert-Schmidt topology is stronger than the weak topology
on $HS$, the subspace $\overline{f^{\,t}(B^{m})}^{w}\cap HS^{n}$ is already
closed in the Hilbert-Schmidt topology.

\n Assume that $\xi \in HS^{n}$ belongs to $(\overline{\ker f}^{HS})^{\perp }$.
Then $\xi \perp \ker f$ and viewing $\xi \in HS^{n}\subset B^{n}$,
as an element of $B^{n}$
we find using the compatibility of the pairings 
that $\xi $ is in the co-kernel of $f$, so that $\xi \in \overline{f^{\,t}(B^{m})}^{w}$
and $\xi \in \overline{f^{\,t}(B^{m})}^{w}\cap HS^{n}$
as claimed. 

\n Finally note that $\overline{f^{\,t}(HS^{m})}^{w}=
\overline{f^{\,t}(B^{m})}^{w}$, because $HS$
is weakly-dense in $B$ and $f^{\,t}$ is weakly continuous.
\end{proof}

\section{First Betti number and $\Delta $.}

\n In this section, we concentrate on the first $L^{2}$-Betti number.

\subsection{$\beta _{1}^{(2)}$ as a limit.\label{sub:beta1inductive}}

\subsubsection{Sub-complexes associated to a set of generators.}

We recall that all $L^{2}$-Betti numbers can be represented as limits,
as described by equation (\ref{eq:betaaslim}). We particularize to
the first Betti number.

\n Let $M$ be a von Neumann algebra with a faithful trace-state $\tau $.

\n Let $F=\{X_{1},\ldots ,X_{n}\}\,,X_{j}\in M$ be a self-adjoint
set of elements in $M$; that is, we assume that $X^* \in F$ 
whenever $X\in F$.  Assume further that $F$ generates $M$ 
as a von Neumann algebra.
Let\[
C_{1}(F)=(M\otimes M^{o})\otimes \Span F\cong (M\otimes M^{o})^{\dim \Span F}\]
and consider\[
\partial _{F}:C_{1}(F)\to M\otimes M^{o},\]
given by
\begin{equation}
\partial_{F} (a\otimes b^{o}\otimes X)=a\,X\otimes b^{o}-a\otimes (X\,b)^{o},\qquad a\otimes b^{o}\in M\otimes M^{o},\  X\in F.
\label{df}\end{equation}
Then\begin{equation}
\ker \partial _{F}\to C_{1}(F)\stackrel{\partial _{F}}{\longrightarrow }M\otimes M^{o}\to M\to 0\label{eq:weakresgenerators}\end{equation}
is a sub-complex of the bar resolution of $M$. The sequence above
is not exact, and the quotient of $M\otimes M^{o}$ by the image
of $\partial _{F}$ is the left $M\otimes M^{o}$-module 
$(M\otimes M^{o})\otimes_{A\otimes A^{o}} A$, where $A$
is the algebra generated by $F$.  When viewed as a bimodule over $M$
this left $M\otimes M^{o}$-module can be identified by the map
$x\otimes y^{o}\mapsto x\otimes y$ with $M\otimes_{A} M$, 
where $(m\otimes n^{o})\cdot( x\otimes_{A}y) = m\,x\otimes_{A}y\,n$. 
We shall thus
use the  notation,
$$
M\otimes_{A} M :=(M\otimes M^{o})\otimes_{A\otimes A^{o}} A
$$
The sequence
\begin{equation}
\ker \partial _{F}\to C_{1}(F)\stackrel{\partial _{F}}{\longrightarrow }M\otimes M^{o}\to M\otimes _{A}M\to 0
\label{res1}\end{equation}
is exact. Applying the induction functor  $(M\bar{\otimes }M^{o})\otimes_{M\otimes M^{o}}$ one gets
 the map\[
1\otimes \partial _{F}:(M\bar{\otimes }M^{o})\otimes \Span F\to M\bar{\otimes }M^{o},\]
given as in (\ref{df}) by right multiplication by $X\otimes 1-1\otimes X^{o}$. Then the first homology of
the induced complex from (\ref{res1}) (or (\ref{eq:weakresgenerators})) is given by\[
H(F)=\frac{\ker (1\otimes \partial _{F})}{(M\bar{\otimes }M^{o})\cdot \ker \partial _{F}}.\]
In other words,

\begin{lem}
\label{lem:HofFasTor}Let $A$ be the algebra generated by $F$. Then
\[
H(F)=\Tor _{1}^{M\otimes M^{o}}(M\otimes _{A}M,M\bar{\otimes }M^{o}).\]
In particular, $H(F)$ depends only on the inclusion $A\subset M$
and not on $F$.
\end{lem}

\n We let $
\beta (F)=\dim _{M\bar{\otimes }M^{o}}H(F)$
so that by Lemma \ref{lem:l2closure} one has
$$
\beta (F) = \beta^{(2)} (\partial _{F})
$$
\n Note that if $F\subset F'$, then there is a natural inclusion map
$i_{F',F}:C_{1}(F)\to C_{1}(F')$. This map induces a homomorphism\[
(i_{F',F})_{*}:H(F)\to H(F').\]

\n For $F\subset F'$ two finite subsets, let\[
H(F:F')=\frac{i_{F',F}(\ker (1\otimes \partial _{F}))}{i_{F',F}(\ker (1\otimes \partial _{F}))\cap (M\bar{\otimes }M^{o})\cdot \ker \partial _{F'}}=(i_{F',F})_{*}H(F).\]
Note that $i_{F',F}(\ker (1\otimes \partial _{F}))\cap (M\bar{\otimes }M^{o})\cdot \ker \partial _{F'}$
is exactly the intersection of $(M\bar{\otimes }M^{o})\otimes \Span F$
with $(M\bar{\otimes }M^{o})\cdot \ker \partial _{F'}$.

\n Then (\ref{eq:betaaslim}) implies that\[
\beta _{1}^{(2)}(M)=\sup _{F}\inf _{F'\supset F}\dim _{M\bar{\otimes }M^{o}}H(F:F').\]
We are thus led to the natural question of the computation of\[
\beta (F:F')=\dim _{M\bar{\otimes }M^{o}}H(F:F')\]
 and, in particular, of $
\beta (F)=\dim _{M\bar{\otimes }M^{o}}H(F)$
(which corresponds to the case that $F=F'$).

\n Note in particular that\[
\beta (F:F')\leq \beta (F)\]
(since the dimension of the image via $(i_{F',F})_{*}$ is not larger
that the dimension of the domain), and also that
\begin{equation}\label{betafprime}
\beta (F:F')\leq \beta (F')
\end{equation}
(since the image via $(i_{F',F})_{*}$ is a sub-module of $H(F')$).

\n Let $F=(X_{1},\ldots ,X_{n})$, $F'=F\cup (Y_{1},\ldots ,Y_{m})$,
and assume for simplicity that $n=\dim \Span F$. Denote by $\partial _{F}^{(2)}$
the extension of $\partial _{F}$ to $(L^{2}(M)\bar{\otimes }L^{2}(M^{o}))^{n}\cong L^{2}(M\bar{\otimes }M^{o})\otimes \Span F$
obtained by continuity.
By Lemma \ref{lem:algsameasFR}, equation (\ref{eq:betaFR}), we have
that\[
\beta (F)=\beta ^{(2)}(\partial _{F})=\dim _{M\bar{\otimes }M^{o}}\ker \partial _{F}^{(2)}-\dim _{M\bar{\otimes }M^{o}}\overline{\ker \partial _{F}^{(2)}\cap FR^{n}}.\]
\medskip

\begin{lem}
\label{lem:kerPartial2}Let $F\subset F'$ be self-adjoint sets 
of elements
in $M$, as before. If $F$ generates $M$ as a von Neumann algebra,
then\[
\dim _{M\bar{\otimes }M^{o}}i_{F',F}\ker \partial _{F}^{(2)}=n-(1-\beta _{0}^{(2)}(M,\tau )),\]
and\begin{eqnarray*}
\beta (F) & = & n-(1-\beta _{0}^{(2)}(M,\tau ))-\dim _{M\bar{\otimes }M^{o}}\overline{(\ker \partial _{F}^{(2)}\cap FR^{n})},\\
\beta (F:F') & = & n-(1-\beta _{0}^{(2)}(M,\tau ))-\dim _{M\bar{\otimes }M^{o}}\overline{\ker \partial _{F'}}\cap i_{F',F}\ker \partial _{F}^{(2)}.
\end{eqnarray*}

\end{lem}
\begin{proof}
We just need to prove the first statement.
The inclusion map
$i_{F',F}:C_{1}(F)\to C_{1}(F')$ is injective, so that we need only to
consider the case that $F=F'$.
As explained in the proof of Proposition \ref{pro:centerbetti0}, 
$\beta _{0}^{(2)}(M,\tau )$ is $1-\dim _{M\bar{\otimes }M^{o}}(\overline{\Image(\partial_1^{(2)})})$ where $\partial_1^{(2)}$ comes from the bar resolution.

\n Considering the kernel and cokernel of the morphism\begin{equation}
 M\bar{\otimes }M^{o}\otimes F\stackrel{\partial _{F}^{(2)}}{\longrightarrow }M\bar{\otimes }M^{o},\label{eq:d1}\end{equation}
it is enough to show that
$$
 \overline{\Image(\partial_{F}^{(2)})}=\overline{\Image(\partial_1^{(2)})}\,.
$$
By construction $\overline{\Image(\partial_1^{(2)})}$
is the strong closure in $M\bar{\otimes }M^{o}$
of the left ideal $L$ generated by \[
V=\{n\otimes 1-1\otimes n^{o}:  n\in M\}.\]

\n We use the map $\Psi :M\bar{\otimes }M^{o}\to HS=HS(L^{2}(M))$ (of
(\ref{eq:defofPsi})); one has by (\ref{bim}) or (\ref{rm})
\begin{equation}\label{psit}
\Psi(x  (n\otimes 1-1\otimes n^{o}))=[Jn^*J,\, \Psi(x)]\,,\quad \forall x\in M\bar{\otimes }M^{o}\,,
n\in M.
\end{equation}
We adopt the following notation for any bounded operator $T$ in $L^{2}(M)$,
\begin{equation}\label{ttt}
\sigma(T):= T^\sigma=JT^*J \,,\quad \forall T\in B\,.
\end{equation}
This gives an antiautomorphism of $B$ which restricts
to the canonical antiisomorphism $\sigma:M\rightarrow M'$.
We thus see that the closure of $L$ in $L^{2}(M\bar{\otimes }M^{o})$
can be identified with the subspace\[
\overline{[HS,M']},\] which is the closure of the 
linear span of commutators of $M'$ with $HS$.

\n Similarly the closure of $\Image(\partial_{F}^{(2)})$ is the subspace 
\[
\overline{[HS,F^\sigma]}\,,\quad F^\sigma:=\sigma(F).\]

\n We just need to show that $[HS,F^\sigma]$ is dense in $[HS,M']$ in the $HS$-topology.
Then the algebra  $A$ generated by $F$ is 
$\ast$-strongly dense in $M$ by hypothesis.
For fixed $T \in HS$ the map $x \mapsto  [T,x]$ is continuous on bounded sets,
from $B$ endowed with the strong topology to $HS$. Thus 
$[HS,\sigma(A)]$ is dense in $[HS,M']$.

\n It remains to sow that $[HS,\sigma(A)]=[HS,F^\sigma]$, which follows from the 
indentities,
\begin{eqnarray*}
[T,X_1 X_2\ldots X_n] & = & [X_2\ldots X_n \,T,\,X_1] +[X_3\ldots X_n \,T\,X_1,\,X_2]\\
 & + & \ldots+[X_{j+1}\ldots X_n \,T\,X_1\ldots X_{j-1},\,X_j]\\
 & + &  \ldots+ [T\,X_1\ldots X_{n-1},\,X_n] 
\end{eqnarray*}
and the fact that $HS$ is a two sided ideal.
\end{proof}

\subsubsection{$\Delta (F)$ and $\Delta (F:F')$.}

It is thus of interest to consider the quantities:\begin{eqnarray}
\Delta (F) & = & n-\dim _{M\bar{\otimes }M^{o}}\overline{(\ker \partial _{F}^{(2)}\cap FR^{n})}=n-\dim _{M\bar{\otimes }M^{o}}\overline{\ker \partial _{F}},\label{eq:deltaF}\\
\Delta (F:F') & = & n-\dim _{M\bar{\otimes }M^{o}}\overline{\ker \partial _{F'}}\cap i_{F',F}\ker \partial _{F}^{(2)},\nonumber \\
\Delta (M,\tau ) & = & \sup _{F\textrm{ s.t. }M=W^{*}(F)}\inf _{F'\supset F}\Delta (F;F'),\nonumber 
\end{eqnarray}
where in the last equation we require that $F$ generates $M$.

\n Explicitly, if $F=(X_{1},\ldots ,X_{n})$, $F'=F\cup (Y_{1},\ldots ,Y_{m})$,
we have:\begin{eqnarray}\label{explicit}
\Delta (F) & = & n-\dim _{M\bar{\otimes }M^{o}}\overline{\{(T_{1},\ldots ,T_{n})\in FR^{n}:\sum [T_{j},X_{j}^\sigma]=0\}}^{HS},\\
\Delta (F:F') & = & n-\dim _{M\bar{\otimes }M^{o}}\Big (HS^{n}\oplus 0\ \nonumber \\
 &  & \cap \  \overline{\{(T_{1},\ldots ,T_{n},S_{1},\ldots ,S_{m})\in FR^{n+m}}\nonumber\\
 &  & \qquad \overline{:\sum [T_{j},X_{j}^\sigma]+\sum [S_{j},Y_{j}^\sigma]=0\}}^{HS}\Big ),\nonumber
\end{eqnarray}
where we used  the map $\Psi :M\bar{\otimes }M^{o}\to HS=HS(L^{2}(M))$.
Note that the $X_j, Y_k$ are moved to the commutant
$M'$ of $M$ by the map $\sigma$, as follows from (\ref{psit}), it is thus clear
that the  subspaces of $HS^n$ involved in the above equations are $M$-bimodules.
In either equation above, $FR^{n}$ can be replaced by $\Psi (M\otimes M^{o})^{n}\subset FR^{n}$.

\n Furthermore, we have by Lemma \ref{lem:kerPartial2}:\begin{equation}
\beta _{1}^{(2)}(M,\tau )=\Delta (M,\tau )-(1-\beta _{0}^{(2)}(M,\tau )).\label{eq:beta1MdeltaM}\end{equation}

\n Note that if $F=(X_{1},\ldots ,X_{n})$, then $\partial _{F}:FR^{n}=FR\otimes \Span F\to FR$
is given by\[
\partial _{F}(T_{1},\ldots ,T_{n})=-\sum [T_{i},X_{i}^\sigma].\]
The transpose of $\partial _{F}$ is the map\[
\partial _{F}^{\,t}:B(L^{2}(M))\to B(L^{2}(M))^{n}\]
given by\begin{equation}
\partial _{F}^{\,t}(D)=([D,X_{1}^\sigma],\ldots ,[D,X_{n}^\sigma]).\label{eq:partialstarF}
\end{equation}
In view of Lemma \ref{lem:dualbetti}, we have the following description
of $\Delta (F)$:\begin{equation}
\Delta (F)=\dim _{M\bar{\otimes }M^{o}}\overline{\partial _{F}^{\,t}(B(L^{2}(M))}^{w}\cap HS^{n}.\label{eq:Deltadualdescription}\end{equation}
Similarly, if $F'=F\cup \{Y_{1},\ldots ,Y_{m}\}$, then\begin{equation}
\Delta (F:F')=\dim _{M\bar{\otimes }M^{o}}\pi _{n}(\overline{\partial _{F'}^{\,t}(B(L^{2}(M))}^{w}\cap HS^{n+m}),\label{eq:dualRelDelta}\end{equation}
where $\pi _{n}:HS^{n+m}\to HS^{n}$ denotes the orthogonal projection
onto the first $n$ coordinates.

\subsection{Properties of $\Delta $.}

We record the following properties of $\Delta $:

\begin{thm}
Let $X_{1},\ldots ,X_{n}\in (M,\tau )$ be a fixed self-adjoint
set of elements.
Then we have:

\n (a) $\Delta (X_{1},\ldots ,X_{n})\leq n$. 

\n (b) $\Delta (X_{1},\ldots ,X_{n})$ depends only on the pair $(A,\tau |_{A})$,
where $A$ is the algebra generated by $X_{1},\ldots ,X_{n}$

\n (c) Let $\Gamma $ be a finitely generated group, and let $X_{1},\ldots ,X_{n}\in L(\Gamma )$
be a family of unitaries associated to a symmetric set of generators of $\Gamma $.
Then\[
\Delta (X_{1},\ldots ,X_{n})\leq \beta _{1}^{(2)}(\Gamma )-\beta _{0}^{(1)}(\Gamma )+1.\]
If in addition $\Gamma $ is combable, we have that\[
\Delta (X_{1},\ldots ,X_{n})=\beta _{1}^{(2)}(\Gamma )-\beta _{0}^{(2)}(\Gamma )+1.\]

\n (d) For all $m<n$,\[
\Delta (X_{1},\ldots ,X_{n})\leq \Delta (X_{1},\ldots ,X_{m})+\Delta (X_{m+1},\ldots ,X_{n}).\]

\n (e) Let $1<m<n$, and assume that the families $X_{1},\ldots ,X_{m}$,
$X_{m+1},\ldots ,X_{n}$ are free. Then\[
\Delta (X_{1},\ldots ,X_{n})=\Delta (X_{1},\ldots ,X_{m})+\Delta (X_{m+1},\ldots ,X_{n}).\]

\end{thm}
\begin{proof}
(a) follows immediately from the definition of $\Delta $.

\n (b) Let $A$ be the algebra generated by $X_{1},\ldots ,X_{n}$, and
let $N$ be the von Neumann algebra generated by $A$ inside of $M$.
Let $F=(X_{1},\ldots ,X_{n})$. By the obvious variant
of Lemma \ref{lem:kerPartial2} for non generating sets, we
find that\[
\Delta (F)=\beta _{1}^{(2)}(F)+(1-\beta _{0}^{(2)}(N)).\]
where $\beta _{1}^{(2)}(F)$ is computed inside $(M,\tau )$.
Moreover, by Lemma \ref{lem:HofFasTor}, we have that\[
\beta _{1}^{(2)}(F)=\dim _{M\bar{\otimes }M^{o}}\Tor _{1}^{M\otimes M^{o}}(M\otimes _{A}M,M\bar{\otimes }M^{o}).\]
Since the functors\[
M\otimes _{N}\cdot ,\qquad \cdot \otimes _{N^{o}}M^{o},\qquad M\bar{\otimes }M^{o}\otimes _{N\bar{\otimes }N^{o}}\]
are flat \cite{luck:hilbertmodules}, it follows that\[
\Tor _{1}^{M\otimes M^{o}}(M\otimes _{A}M,M\bar{\otimes }M^{o})=M\bar{\otimes }M^{o}\otimes _{N\bar{\otimes }N^{o}}\Tor _{1}^{N\otimes N^{o}}(N\otimes _{A}N,N\bar{\otimes }N^{o}).\]
Finally, since\[
\dim _{M\bar{\otimes }M^{o}}(M\bar{\otimes }M^{o}\otimes _{N\bar{\otimes }N^{o}}W)=\dim _{N\bar{\otimes }N^{o}}W,\]
we find that\[
\Delta (F)=\dim _{N\bar{\otimes }N^{o}}\Tor _{1}^{N\otimes N^{o}}(N\otimes _{A}N
,N\bar{\otimes }N^{o})+1-\beta _{0}^{(2)}(N),\]
which depends only on $A$ and $\tau |_{A}$.

\n (c) The inequality follows from Theorem \ref{Prop:sameAsInGroupCase}.
The equality in the combable case follows from Theorem \ref{thm:combableBetti}.

\n (d) Let $F_{1}=(X_{1},\ldots ,X_{m})$ and $F_{2}=(X_{m+1},\ldots ,X_{n})$,
$F=F_{1}\cup F_{2}$. Let $V_{i}=\Span F_{i}$, $C_{i}=(M\bar{\otimes }M^{o})\otimes V_{i}$,
$i=1,2$. Put $C=(M\bar{\otimes }M^{o})\otimes \Span (V_{1},V_{2})$.
Consider\[
\partial _{F_{i}}:C_{i}\to M\otimes M^{o}\]
given by\[
\partial _{F_{i}}(a\otimes b^{o}\otimes x)=a\,x\otimes b^{o}-a\otimes (x\,b)^{o}.\]
Then\[
(\ker \partial _{F_{1}})\oplus (\ker \partial _{F_{2}})\subset \ker \partial _{F}\cap C.\]
Thus\[
\dim _{M\bar{\otimes }M^{o}}\overline{\ker \partial _{F}}\geq \dim _{M\bar{\otimes }M^{o}}\overline{\ker \partial _{F_{1}}}+\dim _{M\bar{\otimes }M^{o}}\overline{\ker \partial _{F_{2}}}.\]
In view of (b), this implies the desired inequality for $\Delta $.

\n (e) Let $M_{1}=W^{*}(X_{1},\ldots ,X_{m})$ and $M_{2}=W^{*}(X_{m+1},\ldots ,X_{n})$.
By Remark 13.2(e) of \cite{dvv:entropy6}, there exist operators $D_{1}$,
$D_{2}$ in $B(L^{2}(M))$ so that 
\begin{equation} \label{dual}
[D_{i},M_{k}]=\{0\},\quad [D_{i},m]=[m,\Psi (1\otimes 1)],
\end{equation}
 for all $i\neq k$ and $m\in M_{i}$. These operators are denoted by
$T_{j}$ in \cite{dvv:entropy6} and are called a dual system to $M_{1},M_{2},\mathbb{C}1$.
It is worth mentioning in conjunction with Def. 13.1 of \cite{dvv:entropy6}
that a single algebra $A$ always has a dual system relative to $\mathbb{C}1$,
namely the operator of orthogonal projection onto $\mathbb{C}1$ in
$L^{2}(A)$.

\n One can in fact explicitly describe these operators. Denote by $H_{i}^{0}$
the space $L^{2}(M_{i})\ominus \mathbb{C}1$. Then since $M=M_{1}*M_{2}$,
\[
L^{2}(M)=\mathbb{C}1\oplus \bigoplus _{k}\bigoplus _{i_{1}\neq i_{2},i_{2}\neq i_{3},\ldots i_{k-1}\neq i_{k}}H_{i_{1}}^{0}\otimes \cdots \otimes H_{i_{k}}^{0}.\]
 We refer the reader to \cite{dvv:book} for more details and the
definition of the action of $M_{j}$ on this space. The operator $D_{k}$
is then given by\[
D_{k}1=1\]
and\[
D_{k}\xi _{1}\otimes \cdots \otimes \xi _{r}=\begin{cases} 0, & \xi _{r}\in H_{k}^{o}\\
 \xi _{1}\otimes \cdots \otimes \xi _{r}, & \textrm{otherwise.}\\ \end{cases}
\]

\n  Let $F=(X_{1},\ldots ,X_{n})$. Assume now that\[
\sum _{i=1}^{n}T_{i}\otimes X_{i}\in \ker \partial _{F},\qquad T_{i}\in M\otimes M^{o}.\]
Then \[
\sum _{i=1}^{n}[\Psi (T_{i}),X_{i}^\sigma]=0.\]
 Let $k$ be equal to $1$ or $2$, and write $I_{1}=\{1,\ldots ,m\}$,
$I_{2}=\{m+1,\ldots ,n\}$. Then for all $a,b\in M$, \begin{eqnarray*}
0 & = & \sum _{i}\Tr ([\Psi (T_{i}),X_{i}^\sigma]\,a\, D_{k}^\sigma\, b)\\
 & = & \sum _{i}\Tr (\Psi (T_{i})[X_{i}^\sigma,\,a\, D_{k}^\sigma\, b])\\
& = & \sum _{i}\Tr (\Psi (T_{i})\,a\,[X_{i}^\sigma,\, D_{k}^\sigma]\, b)\\
 & = & \sum _{i\in I_{k}}\Tr (\Psi (T_{i})[\Psi (a\otimes b^{o}),X_{i}^\sigma])\\
 & = & -\sum _{i\in I_{k}}\Tr ([\Psi (T_{i}),X_{i}^\sigma]\Psi (a\otimes b^{o})).
\end{eqnarray*}
where we used (\ref{dual}) and the equality $\Psi (1\otimes 1)^\sigma=\Psi (1\otimes 1)$.
 
\n It follows (since $a,b\in M$ were arbitrary) that \[
\sum _{i=1}^{m}[\Psi (T_{i}),X_{i}^\sigma]=\sum _{i=m+1}^{n}[\Psi (T_{i}),X_{i}^\sigma]=0.\]

\n It follows that $\sum _{i=1}^{m}T_{i}\otimes X_{i}\in \ker \partial _{F}$, and
$\sum _{i=m+1}^{n}T_{i}\otimes X_{i}\in \ker \partial _{F}$.

\n Let $F_{1}=(X_{1},\ldots ,X_{m})$, $F_{2}=(X_{m+1},\ldots ,X_{n})$.
If we denote by $V_{j}$ the span of $F_{j}$  and  let $V=\Span (V_{1},V_{2})$, then we have shown that
 \[
\ker \partial _{F}\subset \ker \partial _{F}\cap M\otimes M^{o}\otimes V_{1}+\ker \partial _{F}\cap M\otimes M^{o}\otimes V_{2},\]
so that\[
\ker \partial _{F}\subset \ker \partial _{F_{1}}+\ker \partial _{F_{2}}.\]

\n Thus\[
\dim _{M\bar{\otimes }M^{o}}(\overline{\ker \partial _{F}})\leq \sum _{k=1}^{2}\dim _{M\bar{\otimes }M^{o}}(\overline{\ker \partial _{F_{k}}}).\]

\n From this and part (b) we conclude that \[
\Delta (X_{1},\ldots ,X_{n})\geq \Delta (X_{1},\ldots ,X_{m})+\Delta (X_{m+1},\ldots ,X_{n}).\]

\n  Since we have proved the opposite inequality in part (c), the desired
equality now follows. 
\end{proof}

\subsection{$\Delta $ and diffuse center.}

\begin{lem}
\label{lem:diffuseCenter}Let $F=(X,X_{1},\ldots ,X_{n})$, and assume
that $[X,X_{j}]=0$ for all $j$. Assume furthermore that the spectrum
of $X$ is diffuse. Then $\Delta (F)=1$.
\end{lem}
\begin{proof}
Let $M=W^{*}(F)$. Since $\Delta (F)=\beta _{1}^{(2)}+(1-\beta _{0}^{(2)}(M))$,
we find that\[
\Delta (F)\geq 1-\beta _{0}^{(2)}(M).\]
Since $M$ contains a diffuse von Neumann algebra (namely, $W^{*}(X)$),
it follows from Proposition \ref{pro:centerbetti0} that $\beta _{0}^{(2)}(M)=0$.
Thus $\Delta (F)\geq 1.$

\n  For the opposite inequality, let $FR$ be the ideal of finite-rank
operators on $L^{2}(M)$, and assume that $Q_{1},\ldots ,Q_{n}\in FR^{n}$
are arbitrary. Let\begin{eqnarray*}
T_{j} & = & [Q_{j},X^\sigma],\\
T & = & -\sum _{j=1}^{n}[Q_{j},X_{j}^\sigma].
\end{eqnarray*}
Then since $[X_{j}^\sigma,X^\sigma]=0$ for all $j$, we have that\begin{eqnarray*}
[T,X^\sigma]+\sum _{j=1}^{n}[T_{j},X_{j}^\sigma] & = & 
\sum _{j=1}^{n}-[[Q_{j},X_{j}^\sigma],X^\sigma]+[[Q_{j},X^\sigma],X_{j}^\sigma]=0\,\\
\end{eqnarray*}
by the Jacobi identity. Thus the image of the map\[
FR^{n}\ni (Q_{1},\ldots ,Q_{n})\mapsto (T,T_{1},\ldots ,T_{n})\]
lies inside $\ker \partial _{F}$ (we identify as usual $M\otimes M^{o}$
with a subset of $FR$ via the map $\Psi $; see Lemma \ref{lem:algsameasFR}).
It follows that the closure of $\ker \partial _{F}$ in the Hilbert-Schmidt
norm contains the image of the map\[
\phi :HS^{n}\ni (Q_{1},\ldots ,Q_{n})\mapsto (-\sum _{j=1}^{n}[Q_{j},X_{j}^\sigma],[Q_{1},X^\sigma],\ldots ,[Q_{n},X^\sigma]).\]
Since $X$ has diffuse spectrum, the commutant of $W^{*}(X^\sigma)$ in $B(L^{2}(M))$
does not intersect compact (and hence Hilbert-Schmidt) operators.
Thus the map $\phi $ is injective. Hence by Luck's results on additivity
of dimension for weakly exact sequences \cite{luck:foundations1}
we conclude that\[
\dim _{M\bar{\otimes }M^{o}}\overline{\ker \partial _{F}}^{HS}\geq \dim _{M\bar{\otimes }M^{o}}\Image \phi =\dim _{M\bar{\otimes }M^{o}}HS^{n}=n.\]
Thus\[
\Delta (F)=n+1-\dim _{M\bar{\otimes }M^{o}}\overline{\ker \partial _{F}}^{HS}\leq 1.\]
We conclude that $\Delta (F)=1$.
\end{proof}
\begin{cor}
If $(M,\tau )$ is a von Neumann algebra, and $M$ has a diffuse center,
then $\beta _{1}^{(2)}(M,\tau )=0$ and $\Delta (M,\tau )=1$.
\end{cor}
\begin{proof}
Let $X$ be a generator of the center of $M$. Then for any finite
subset $F$ of self-adjoint elements of $M$, we have that\[
\Delta (F\cup \{X\})=1,\]
by Lemma \ref{lem:diffuseCenter}. Hence if $F$ generates $M$, $F'\supset F$
and $X\in F'$, then by (\ref{betafprime})
\begin{eqnarray*}
\Delta (F:F') & = & \beta _{1}^{(2)}(F:F')+(1-\beta _{0}^{(2)}(M))\\
 & \leq  & \beta _{1}^{(2)}(F')+(1-\beta _{0}^{(2)}(M))\\
 & = & \Delta (F')=1.
\end{eqnarray*}

\n  Thus for any $F$ generating $M$, we have that\[
\inf _{F'\supset F}\Delta (F:F')\leq \Delta (F:F\cup \{X\})=1,\]
so that $\Delta (M,\tau )\leq 1$. Since $M$ contains a diffuse von
Neumann algebra (namely, $W^{*}(X)$), it follows that $\Delta (M,\tau )=1$,
and that $\beta _{1}^{(2)}(M)=0$.
\end{proof}

\section{$\Delta $ and free entropy dimension.}

\subsection{Free entropy dimension.}

The properties of $\Delta (X_{1},\ldots ,X_{n})$ seem very similar
to those enjoyed by the various versions of Voiculescu's free entropy
dimension. We therefore are interested in connections between the
two quantities.

\subsubsection{Non-microstates entropy dimension.}

We consider the free entropy dimensions defined in terms of the non-microstates
free entropy and the non-microstates free Fisher information. Let
$S_{1},\ldots ,S_{n}$be a free semicircular family, free from $(X_{1},\ldots ,X_{n})$.
Then set\[
\delta ^{*}(X_{1},\ldots ,X_{n})=n-\liminf _{\varepsilon \to 0}\frac{\chi ^{*}(X_{1}+\sqrt{\varepsilon }S_{1},\ldots ,X_{n}+\sqrt{\varepsilon }S_{n})}{\log \varepsilon ^{1/2}}\]
and\[
\delta ^{\star }(X_{1},\ldots ,X_{n})=n-\liminf _{\varepsilon \to 0}\varepsilon\, \Phi ^{*}(X_{1}+\sqrt{\varepsilon }S_{1},\ldots ,X_{n}+\sqrt{\varepsilon }S_{n}).\]
Note that $\delta ^{\star }$ is obtained from $\delta ^{*}$ by formally
applying L'Hopital's rule to the limit. We will also use the microstates
free entropy dimension $\delta $ and $\delta _{0}$, which were introduced
in \cite{dvv:entropy2,dvv:entropy3}. Here\[
\Phi ^{*}(X_{1}+\sqrt{\varepsilon }S_{1},\ldots ,X_{n}+\sqrt{\varepsilon }S_{n})=\sum _{i=1}^{n}\Vert \xi _{i}^{\varepsilon }\Vert _{2}^{2},\]
where $\xi _{i}^{\varepsilon }\in L^{2}(W^{*}(X_{1},\ldots ,X_{n}))$
are the conjugate variables\[
\xi _{i}^{\varepsilon }=J(X_{i}+\sqrt{\varepsilon }S_{i}:X_{1}+\sqrt{\varepsilon }S_{1},\ldots ,\widehat{X_{i}+\sqrt{\varepsilon }S_{i}},\ldots ,X_{n}+\sqrt{\varepsilon }S_{n})\]
(here $\hat{\cdot }$ denotes omission). 

\n Let\[
E_{\varepsilon }=E_{W^{*}(X_{1}+\sqrt{\varepsilon }S_{1},\ldots ,X_{n}+\sqrt{\varepsilon }S_{n})}\]
be the unique conditional expectation. By \cite{dvv:entropy5}, one
has:\[
\xi _{i}^{\varepsilon }=\frac{1}{\sqrt{\varepsilon }}E_{\varepsilon }(S_{i}).\]
Thus\[
\delta ^{\star }=n-\liminf _{\varepsilon \to 0}\sum _{i=1}^{n}\Vert E_{\varepsilon }(S_{i})\Vert _{2}^{2}.\]

\begin{lem}
\label{lemma:StarvsStar}$\delta ^{\star }(X_{1},\ldots ,X_{n})\geq \delta ^{*}(X_{1},\ldots ,X_{n})$.
\end{lem}
\begin{proof}
Assume that $\delta ^{\star }(X_{1},\ldots ,X_{n})< n-C$. Thus\[
\liminf _{\varepsilon \to 0}\varepsilon \,\Phi ^{*}(X_{1}+\sqrt{\varepsilon }S_{1},\ldots ,X_{n}+\sqrt{\varepsilon }S_{n})> C.\]
Then for some $\varepsilon _{0}>0$ and all $0<\varepsilon <\varepsilon _{0}$,
we have that\[
\varepsilon \,\Phi ^{*}(X_{1}+\sqrt{\varepsilon }S_{1},\ldots ,X_{n}+\sqrt{\varepsilon }S_{n})\geq C,\]
so that\[
\Phi ^{*}(X_{1}+\sqrt{\varepsilon }S_{1},\ldots ,X_{n}+\sqrt{\varepsilon }S_{n})\geq \frac{C}{\varepsilon }.\]
Thus for all $0<\varepsilon <\varepsilon _{0}$, \begin{eqnarray*}
\frac{1}{2}\int _{\varepsilon }^{\varepsilon _{0}}\Phi ^{*}(X_{1}+\sqrt{t}S_{1},\ldots ,X_{n}+\sqrt{t}S_{n})\,dt  & \geq  & \frac{1}{2}\int _{\varepsilon }^{\varepsilon _{0}}\frac{C}{\  t}dt\\
 & = & C(\log \varepsilon _{0}^{1/2}-\log \varepsilon ^{1/2}).
\end{eqnarray*}
Now $\chi _{\varepsilon }=\chi ^{*}(X_{1}+\sqrt{\varepsilon }S_{1},\ldots ,X_{n}+\sqrt{\varepsilon }S_{n})$
is given by (cf \cite{dvv:entropy5})\begin{eqnarray*}
\chi _{\varepsilon } & =\frac{1}{2} & \int _{0}^{\infty }\left(\frac{n}{1+t}-\Phi ^{*}(X_{1}+\sqrt{t+\varepsilon }\,S_{1},\ldots ,X_{n}+\sqrt{t+\varepsilon }\,S_{n}\right)dt\\
 & = & \frac{1}{2}\int _{\varepsilon }^{\infty }\left(\frac{n}{1+t-\varepsilon }-\Phi ^{*}(X_{1}+\sqrt{t}S_{1},\ldots ,X_{n}+\sqrt{t}S_{n}\right)dt\\
 & \leq  & K+\frac{1}{2}\int _{\varepsilon }^{\varepsilon _{0} }\left(\frac{n}{1+t-\varepsilon }-\Phi ^{*}(X_{1}+\sqrt{t}S_{1},\ldots ,X_{n}+\sqrt{t}S_{n}\right)dt,
\end{eqnarray*}
for some constant $K$ depending only on $\varepsilon _{0}$ and $X_{1},\ldots ,X_{n}\,.$

\n Thus 
$$
\chi _{\varepsilon }\leq K+\frac{n}{2}\log \left(1+\varepsilon _{0}-\varepsilon \right)+C(\log \varepsilon ^{1/2}-\log \varepsilon _{0}^{1/2})
$$
 Since for small $\varepsilon $, $\log \varepsilon $
is negative, it follows that\[
\liminf _{\varepsilon \to 0}\frac{\chi ^{*}(X_{1}+\sqrt{\varepsilon }S_{1},\ldots ,X_{n}+\sqrt{\varepsilon }S_{n})}{\log \varepsilon ^{1/2}}\geq \liminf _{\varepsilon \to 0}\frac{C\log \varepsilon ^{1/2}}{\log \varepsilon ^{1/2}}=C.\]
Thus $\delta ^{*}(X_{1},\ldots ,X_{n})\leq n-C$. Since $C$ is arbitrary,
we get that $\delta ^{*}(X_{1},\ldots ,X_{n})\leq \delta ^{\star }(X_{1},\ldots ,X_{n})$.
\end{proof}

\subsection{The inequality $\Delta \geq \delta ^{\star }$.}

In preparation for the next theorem, we need to set up some notation.

\n Let $S_{1},\ldots ,S_{n}$ be a free semicircular system, free from
$(X_{1},\ldots ,X_{n})$. Let $X_{j}({\varepsilon })=X_{j}+\sqrt{\varepsilon }\,S_{j}$,
$j=1,\ldots ,n$. Let also $M_{\varepsilon }=W^{*}(X_{1}({\varepsilon }),\ldots ,X_{n}({\varepsilon }))$,
$N=W^{*}(X_{1},\ldots ,X_{n},S_{1},\ldots ,S_{n})$, $H=L^{2}(N)$.
(We recall in Appendix II  the details of this standard construction  (cf. \cite{dvv:entropy5}))
Thus $M_{\varepsilon }\subset N\subset B(H)$. Let $E_{\varepsilon }$
be the orthogonal projection from $H$ onto $L^{2}(M_{\varepsilon })$.
We denote by the same symbol the conditional expectation from $N$
onto $M_{\varepsilon }$. Note that $P_{1}=E_{\varepsilon }P_{1}=P_{1}E_{\varepsilon }$.

\n Let $D_{j}({\varepsilon })$, $j=1,\ldots ,n$, be a dual system to
$(X_{1}({\varepsilon }),\ldots ,X_{n}({\varepsilon }))$ on $B(H)$.
That is, we require that 
\begin{equation}\label{dual2}
D_{j}({\varepsilon })=E_{\varepsilon }D_{j}({\varepsilon })E_{\varepsilon }\,,\quad
[D_{j}({\varepsilon }),X_{i}({\varepsilon })]=\delta _{ij}P_{1}
\end{equation}
Such a system is always possible to find: one can set 
\begin{equation}\label{dual3}
D_{j}({\varepsilon })=E_{\varepsilon }\frac{1}{\sqrt{\varepsilon }}Q_{j}E_{\varepsilon }
\end{equation}
where $Q_{j}$ is a right creation operator (cf. Appendix II).
Note that one has the property that $\Vert D_{j}({\varepsilon })\Vert _{\infty }\leq 1/\sqrt{\varepsilon }$. 

\n Identify now $M\otimes M^{o}$ via the map $\Psi $ with a subspace of the space
of finite-rank operators on $H$. Let $T_{1},\ldots ,T_{n}\in \Psi (M\otimes M^{o})$
be so that $\sum [T_{i},X_{i}^{\sigma}(\varepsilon)]=0$. To be explicit,
let $T_{i}=\sum _{k}a_{k}^{i}\,P_{1}\,b_{k}^{i}$, $a_{k}^{i},b_{k}^{i}\in M$. 

\n We first need a few lemmas.

\begin{lem}
\label{Lemma:TandTprime}Let $\delta >0$. Then there exist $x_{k}^{i}(\varepsilon),
y_{k}^{i}(\varepsilon)\in M_{\varepsilon }$ and $\varepsilon _{0}>0$
such that,
\[
\Vert T_{i}'(\varepsilon)-T_{i}\Vert _{HS}\leq \Vert T_{i}'(\varepsilon)-T_{i}\Vert _{1}<\delta \,,\quad \forall \varepsilon <\varepsilon _{0} \]
where \[
T_{i}'(\varepsilon)=\sum _{k}\,x_{k}^{i}(\varepsilon)\,P_{1}\,y_{k}^{i}(\varepsilon)\,.\]

\end{lem}
\begin{proof}
(of Lemma). It is sufficient to prove the statement for single rank-one
operator $T_{a,b}=\,a\,P_{1}\,b\,$. Note that $\Vert T_{a,b}\Vert _{1}=\sup _{\Vert S\Vert _{\infty }=1}|\langle a,Sb\rangle |=\Vert a\Vert _{2}\Vert b\Vert _{2}$. We can assume that 
$\Vert a\Vert _{2}=\Vert b\Vert _{2}=1$.
Choose non-commutative polynomials $p$ and $q$ so that\[
\Vert p(X_{1},\ldots ,X_{n})-a\Vert _{2}\leq \delta /4,\  \Vert q(X_{1},\ldots ,X_{n})-b\Vert _{2}<\delta /4.\]
 Let $\varepsilon _{0}>0$ so that whenever $X_{j}'\in N$
and $\Vert X_{j}-X_{j}'\Vert _{\infty }<2\sqrt{\varepsilon _{0}}$, we have\[
\Vert p(X_{1},\ldots ,X_{n})-p(X_{1}',\ldots ,X_{n}')\Vert _{\infty }<\delta /4,\  \Vert q(X_{1},\ldots ,X_{n})-q(X_{1}',\ldots ,X_{n}')\Vert _{\infty }<\delta /4.\]
Set $x(\varepsilon)=p(X_{1}(\varepsilon),\ldots ,X_{n}(\varepsilon))$,
$y(\varepsilon)=q(X_{1}(\varepsilon),\ldots ,X_{n}(\varepsilon))$.
Let $0<\varepsilon <\varepsilon _{0}$. Then $\,\Vert \,x(\varepsilon)-a\Vert _{2}\leq \delta /2$
and $\Vert y(\varepsilon)-b\Vert _{2}\leq \delta /2$ which gives the answer.
\end{proof}

\medskip
The following lemma is implicit in \cite{dvv:entropy5}, but we restate it
for convenience.
\begin{lem}
\label{lemma:DvsJ}Let $N$ be von Neumann algebra, and let $\tau $
be a faithful normal trace on $N$. Let $H=L^{2}(N,\tau )$, and let $J$
be the Tomita conjugation associated to $N$. Denote by $P_{1}$ the
orthogonal projection onto $1\in H$.

\n Let $Q\in B(H)$ and  $Z\in N$. Then\[
\Tr (P_{1}[Q,Z^{*}])=\langle (JQ^{*}J-Q)1,Z1\rangle .\]

\end{lem}
\begin{proof}
Using $\langle J \,\xi,\,\eta \rangle=\langle J \,\eta,\,\xi \rangle$ we get:
\begin{eqnarray*}
\Tr (P_{1}[Q,Z^{*}]) & = & \langle QZ^{*}1,1\rangle -\langle Z^{*}Q1,1\rangle \\
 & = & \langle Z^{*},Q^{*}1\rangle -\langle Q1,Z\rangle \\
 & = & \langle JQ^{*}J1,Z\rangle -\langle Q1,Z\rangle \\
 & = & \langle (JQ^{*}J-Q)1,Z\rangle ,
\end{eqnarray*}
which is the desired identity.
\end{proof}
\begin{thm}
\label{thm:DeltaGEQdeltastar}$\Delta (X_{1},\ldots ,X_{n})\geq \delta ^{\star }(X_{1},\ldots ,X_{n}).$
\end{thm}
\begin{proof}
Let $T_{1}^{j},\ldots ,T_{n}^{j}$, $j=1,\ldots ,n$ be in $\Psi (M\otimes M^{o})$
and such that $\sum _{i}[T_{i}^{j},X_{i}^\sigma]=0$ for all $j$. Let
$$T_{i}^{j}=\sum _{k}\,a_{k}^{ij}\,P_{1}\,b_{k}^{ij}\,.$$ Then
using (\ref{dual2}), $\sum _{i}[T_{i}^{j},X_{i}^\sigma]=0$, and $\sigma(P_{1})=P_{1}$, we get,
\begin{eqnarray*}
\sum _{j}\langle T_{j}^{j},P_{1}\rangle  & = & \sum _{ij}\Tr (T_{i}^{j}\,[X_{i}^\sigma(\varepsilon),D^\sigma_{j}(\varepsilon)])\\
 & = & \sum _{ij}\Tr ([T_{i}^{j},\,X_{i}^\sigma(\varepsilon)]D^\sigma_{j}(\varepsilon))\\
 & = & \sum _{ij}\Tr ([T_{i}^{j},\,\sqrt{\varepsilon }\,S_{i}^\sigma]\,D^\sigma_{j}(\varepsilon))\\
 & = & \sum _{ij}\Tr (T_{i}^{j}[S_{i}^\sigma,\,\sqrt{\varepsilon }\,D^\sigma_{j}(\varepsilon)]).
\end{eqnarray*}
Now, fix $\kappa >0$ and let $\delta =\kappa /8n^{2}$. Choose $\varepsilon _{0}$
and $T_{ij}'(\varepsilon)$ as in Lemma \ref{Lemma:TandTprime}, so that 
\begin{equation}\label{approx}
\Vert T_{i}^{j}-T_{ij}'(\varepsilon)\Vert _{1}<\delta .
\end{equation}
Then for all $\varepsilon <\varepsilon _{0}$, since $\Vert [\sqrt{\varepsilon }D_{j}(\varepsilon),S_{i}]\Vert _{\infty }\leq 4$,
we find that
\begin{equation}\label{inf1}
\vert\sum _{j}\langle T_{j}^{j},P_{1}\rangle \vert  \leq   \kappa /2+\vert\sum _{ij}\Tr (T_{ij}'(\varepsilon)[\sqrt{\varepsilon }\,D_{j}^\sigma(\varepsilon),S_{i}^\sigma])\vert .
\end{equation}
Since $$T_{ij}'(\varepsilon)=\sum \,x_{k}^{ij}(\varepsilon)\,P_{1}\,y_{k}^{ij}(\varepsilon)
\,,\qquad  x_{k}^{ij}(\varepsilon),y_{k}^{ij}(\varepsilon)\in M_{\varepsilon }$$ and
$P_{1}=E_{\varepsilon }P_{1}E_{\varepsilon }$, we have that $$T_{ij}'(\varepsilon)
=E_{\varepsilon }T_{ij}'(\varepsilon)E_{\varepsilon }.$$
Let $\xi _{i}(\varepsilon)=E_{\varepsilon }(\frac{1}{\sqrt{\varepsilon }}S_{i})$ (cf.\cite{dvv:entropy5})
then $\xi _{i}(\varepsilon)\in M_{\varepsilon }$ 
and $\Vert \sqrt{\varepsilon }\,\xi _{i}(\varepsilon)\Vert _{\infty }\leq \Vert S_{i}\Vert _{\infty }\leq 2$, in fact
$$
E_{\varepsilon }S_{i} E_{\varepsilon }= \sqrt{\varepsilon }\,\xi _{i}(\varepsilon)
$$
Note that $E_{\varepsilon }^\sigma=E_{\varepsilon }$ and $\sqrt{\varepsilon }\,D_{j}^\sigma(\varepsilon)=E_{\varepsilon } Q_{j}^\sigma E_{\varepsilon }$
by (\ref{dual3}). One has\begin{eqnarray}\label{inf2}
 \sum _{ij}\Tr (T_{ij}'(\varepsilon)[\sqrt{\varepsilon }\,D_{j}^\sigma(\varepsilon),S_{i}^\sigma])& = & \sum _{ij}\Tr (T_{ij}'(\varepsilon)[Q_{j}^\sigma,E_{\varepsilon }S_{i}^\sigma E_{\varepsilon }])\\
 & = & \sum _{ij}\Tr (T_{ij}'(\varepsilon)[Q_{j}^\sigma,\sqrt{\varepsilon }\,\xi _{i}(\varepsilon)^\sigma
])\nonumber 
\end{eqnarray}
We thus get, using (\ref{approx}), the fact that $\,a_{k}^{ij}\,$, $\,b_{k}^{ij}\,$
 commute with $\xi _{i}(\varepsilon)^\sigma$
and the inequality $\Vert \sqrt{\varepsilon }\, \xi _{i}(\varepsilon)\Vert _{\infty }\leq \Vert S_{i}\Vert _{\infty }\leq 2$:\begin{eqnarray*}
\vert\sum _{ij}\Tr (T_{ij}'(\varepsilon)[Q_{j}^\sigma,\sqrt{\varepsilon }\,\xi _{i}(\varepsilon)^\sigma
])\vert
 & \leq  &  \kappa/2 +|\sum _{ij}\Tr (T_{i}^{j}
[Q_{j}^\sigma,\sqrt{\varepsilon }\,\xi _{i}(\varepsilon)^\sigma
])\\
 & = & \kappa/2 +|\sum _{ij}\Tr (P_{1}\  [\sum _{k}\,b_{k}^{ij}\,Q_{j}^\sigma\,a_{k}^{ij}\,,\  \sqrt{\varepsilon }\,\xi _{i}(\varepsilon)^\sigma
])|\\
 & = & \kappa/2 +|\sum _{ij}\langle (Y_{ij}-JY_{ij}^{*}J)1,
\sqrt{\varepsilon }\,\xi _{i}(\varepsilon)^\sigma
\rangle _{H}|,
\end{eqnarray*}
 where  $Y_{ij}=\sum _{k}\,b_{k}^{ij}\,Q_{j}^\sigma\,a_{k}^{ij}\,$,
and the last equality is by Lemma \ref{lemma:DvsJ}. 

\n Combining this with (\ref{inf1}) and (\ref{inf2}) we get
\begin{equation}\label{inf3}
\vert\sum _{j}\langle T_{j}^{j},P_{1}\rangle \vert  \leq   \kappa +|\sum _{ij}\langle (Y_{ij}-JY_{ij}^{*}J)1,
\sqrt{\varepsilon }\,\xi _{i}(\varepsilon)^\sigma
\rangle _{H}|.
\end{equation}

\n Let $\eta _{ij}=(Y_{ij}-JY_{ij}^{*}J)1 \in H$. Computing explicitly, we
get that (cf. Appendix II)
\begin{equation}\label{comp}
\eta _{ij}=-\sum _{k}b_{k}^{ij}S_{j}a_{k}^{ij}
\end{equation}
and in particular $\Vert \sum _{j}\eta _{ij}\Vert _{2}^{2}=\sum _{j}\Vert T_{i}^{j}\Vert _{2}^{2}$,
since the subspaces $MS_{j}M$ are orthogonal for $j=1,\ldots ,n$,
and the map $\sum x_{k}\otimes y_{k}\to \sum x_{k}S_{j}y_{k}$ is
an isometry from $L^{2}(M)\bar{\otimes }L^{2}(M)$ into $L^{2}(N)$,
for each $j$ (cf. Appendix II). We thus conclude that\begin{eqnarray*}
|\sum _{j}\langle T_{j}^{j},P_{1}\rangle | & \leq  & \kappa +|\sum _{i}\langle \sum _{j}\eta _{ij},\sqrt{\varepsilon }\,\xi _{i}(\varepsilon)^\sigma
\rangle |\\
 & \leq  & \kappa +(\sum _{i}\Vert \sum _{j}\eta _{ij}\Vert ^{2})^{1/2}(\varepsilon \sum _{i}\Vert \xi _{i}(\varepsilon)\Vert ^{2})^{1/2}\\
 & = & \kappa +(\sum _{ij}\Vert T_{i}^{j}\Vert _{2}^{2})^{1/2}(\varepsilon \Phi ^{*}(X_{1}(\varepsilon),\ldots ,X_{n}(\varepsilon)))^{1/2}.
\end{eqnarray*}
since the free Fisher information is defined as $\Phi ^{*}(X_{1}(\varepsilon),\ldots ,X_{n}(\varepsilon))=\sum _{i=1}^{n}\Vert \xi _{i}(\varepsilon)\Vert _{2}^{2}$.
Passing to $\liminf _{\varepsilon \to 0}$ and noticing that $\kappa $
is arbitrary finally gives us:
\begin{equation}\label{est}
|\sum _{j}\langle T_{j}^{j},P_{1}\rangle |\leq (n-\delta ^{\star }(X_{1},\ldots ,X_{n}))^{1/2}\left(\sum _{ij}\Vert T_{i}^{j}\Vert _{2}^{2}\right)^{1/2}.
\end{equation}

\n The conclusion of the proof of the theorem now follows from the next lemma 
(cf. \cite[Lemma 2.9]{shlyakht:qdim})
applied to the 
von Neumann algebra $M\bar{\otimes }M^{o}$ and the subspace $K$
of  $L^{2}(M\bar{\otimes }M^{o})^n=HS^n$
 closure of the space 
$\{(T_{1},\ldots ,T_{n})\in \Psi (M\otimes M^{o})^{n}:\sum [T_{i},X_{i}^\sigma]=0\}$.\end{proof}

\begin{lem}
\label{Lemma:distAndDim}Let $N$ be a finite von Neumann algebra
with a faithful normal trace $\tau $. Let $n$ be a finite integer,
and let $H=L^{2}(N,\tau )^{n}$ viewed as  a left module over $N$. Denote
by $\Omega \in L^{2}(N,\tau )$ the GNS vector associated to $\tau $. 

\n Let $K\subset H$ be a closed $N$-invariant subspace of $H$. Endow
$M_{n\times n}(L^{2}(N))$ with the norm\[
\Vert h\Vert _{M_{n}}^{2}=\sum _{ij=1}^{n}\Vert h_{ij}\Vert ^{2}.\]
Let $A(K)=\{T\in M_{n\times n}(N):TH\subset K\}\cong K^{n}$. Then
we have:\[
\dim _{N}K={\rm Sup}\vert\langle T, I\rangle \vert^2 /\Vert T \Vert^2\,,\quad T \in A(K),\]
where $I\in M_{n}(H)$ denotes the matrix $I_{ij}=\delta _{ij}\Omega $
\end{lem}
\begin{proof} (of Lemma).
We identify the commutant $N'$ of $N$ acting on $H$  with
the algebra of $n\times n$ matrices $M_{n}(N)$. Endow this algebra
with the non-normalized trace $\Tr $, defined by the property that
$\Tr (1)=n$, where $1\in M_{n}(N)$ denote the identity matrix. Let
$e_{K}\in N'$ be the orthogonal projection from $H$ onto $K$. Then\[
\dim _{N}K=\Tr (e_{K}).\]
Now, $L^{2}(M_{n}(N),\Tr )=M_{n}(H)$ isometrically. Moreover, 
the orthogonal projection of $I$ onto $A(K)$ is $e_{K}\in A(K)$, since $1-e_{K}$
is orthogonal to $A(K)=e_{K}\,M_{n}(N)$. The above supremum is thus reached for $T=e_{K}$
and its  value is $\Tr(e_{K})$ which gives the result.\end{proof}

\subsubsection{Some consequences for $\Delta $.}

\begin{cor}
We have\[
\Delta (X_{1},\ldots ,X_{n})\geq \delta ^{\star }(X_{1},\ldots ,X_{n})\geq \delta ^{*}(X_{1},\ldots ,X_{n})\geq \delta (X_{1},\ldots ,X_{n})\geq \delta _{0}(X_{1},\ldots ,X_{n}).\]

\end{cor}
\n This is immediate from the preceding discussion and the work of Biane,
Guionnet and Capitaine \cite{guionnet-biane-capitaine:largedeviations}.

\n The following corollary gives a strong indication that the first $L^{2}$-Betti
number of a free group factor does not vanish (compare with equations
(\ref{eq:beta1MdeltaM}) and (\ref{eq:deltaF})). 

\begin{cor}
Let $F=(X_{1},\ldots ,X_{n})$ be a self-adjoint
  finite subset of $M$, and assume
that $F$ generates
$M$. 

\n Assume that the microstates free entropy $\chi (X_{1},\ldots ,X_{n})$
is finite. 

\n Then for any seld-adjoint subset $F'$ of $M$, we
have\[
\Delta (F\cup F')\geq n.\]

\end{cor}
\begin{proof}
Let $F'=(Y_{1},\ldots ,Y_{n})$. Then\[
\Delta (F\cup F')\geq \delta (X_{1},\ldots ,X_{n},Y_{1},\ldots ,Y_{m})\geq n\]
where the second inequality follows from \cite{dvv:entropy2}.
\end{proof}
It is of course of interest if one has $\Delta =\delta ^{\star }$.
In conjunction with this, we note the following. Let $F=(X_{1},\ldots ,X_{n})$
be a finite self-adjoint subset of $M$. Consider
as in (\ref{eq:partialstarF})\[
\partial _{F}^{\,t}:B(L^{2}(M))\to B(L^{2}(M))^{n}\]
given by\[
\partial _{F}^{\,t}(D)=([D,X_{1}^\sigma],\ldots ,[D,X_{n}^\sigma]).\]

\begin{thm}
One has\[
\dim _{M\bar{\otimes }M^{o}}\overline{\partial _{F}^{\,t}(B(L^{2}(M))\cap HS^{n}}^{HS}\leq \delta ^{*}(X_{1},\ldots ,X_{n})\]
\begin{eqnarray*}
\delta ^{\star }(X_{1},\ldots ,X_{n}) & \leq  & \dim _{M\bar{\otimes }M^{o}}\overline{\partial _{F}^{\,t}(B(L^{2}(M))}^{w}\cap HS^{n}=\Delta (F).
\end{eqnarray*}

\end{thm}
\begin{proof}
The first inequality is the statement of Corollary 2.12 in \cite{shlyakht:qdim}.
The second inequality is the statement of Theorem \ref{thm:DeltaGEQdeltastar},
together with the {}``dual'' description of $\Delta (F)$ given
in equation (\ref{eq:Deltadualdescription}).
\end{proof}

\subsubsection{Some consequences for free entropy dimension.}

Let $C(\Gamma)$ be the cost of a discrete group $\Gamma$ in the sense of 
\cite{gaboriau:cost}.

\begin{cor}
Let $\Gamma $ be a finitely generated group with a symmetric set of generators $\gamma _{1},\ldots ,\gamma _{n}$.
Denote by $u_{i}=\lambda (\gamma _{i})\in L(\Gamma )$ the corresponding
unitaries in the left regular representation. Let $X_{i}=u_{i}+u_{i}^{*}$,
$Y_{i}=i(u_{i}-u_{i}^{*})$. Then\[
\delta ^{*}(X_{1},\ldots ,X_{n},Y_{1},\ldots ,Y_{n})\leq \beta _{1}^{(2)}(\Gamma )-\beta _{0}^{(2)}(\Gamma )+1 \leq C(\Gamma).\]

\end{cor}
\begin{proof}
The inequality between the cost and $\beta_1 - beta_0 + 1$ is due to Gaboriau
\cite{gaboriau:ell2}.  
The rest of the inequalities follow immediately from 
the corresponding estimate for $\Delta $.  
\end{proof}
\begin{cor}
Let $\Gamma $ be a discrete group with Kazhdan's property (T). Let
$\gamma _{1},\ldots ,\gamma _{n}$ be a symmetric set of generators of $\Gamma $, and
let $u_{j}=\lambda (\gamma _{j})\in L(\Gamma )$ be the associated
unitaries in the left regular representation. Let $X_{i}=u_{i}+u_{i}^{*}$,
$Y_{i}=i(u_{i}-u_{i}^{*})$. Then\[
\delta _{0}(\Gamma )\leq \delta (X_{1},\ldots ,X_{n},Y_{1},\ldots ,Y_{n})\leq \delta ^{\star }(X_{1},\ldots ,X_{n},Y_{1},\ldots ,Y_{n})\leq 1.\]
If moreover $L(\Gamma )$ is diffuse, one has\[
\delta ^{*}(X_{1},\ldots ,X_{n},Y_{1},\ldots ,Y_{n})=\delta ^{\star }(X_{1},\ldots ,X_{n},Y_{1},\ldots ,Y_{n})=1.\]
If $L(\Gamma )$ is diffuse and moreover $L(\Gamma )$ can be embedded
into the ultrapower of the hyperfinite II$_{1}$-factor, one has\[
\delta _{0}(\Gamma )=1.\]

\end{cor}
\begin{proof}
The upper estimates are a consequence of the fact that if $\Gamma $
has property (T), then $\beta _{1}^{(2)}(\Gamma )=0$ (see e.g. \cite{cheeger-gromov:l2,bekka-valette:l2cohomology}).
Thus\[
\Delta (\Gamma )\leq \beta _{1}^{(2)}(\Gamma )-\beta _{0}^{(2)}(\Gamma )+1=1-\beta _{0}^{(2)}(\Gamma )\leq 1.\]
The lower estimate for $\delta ^{*}$ is a consequence of \cite[Theorem 2.13]{shlyakht:qdim}.
The corresponding estimate for $\delta _{0}$ is a consequence of
hyperfinite monotonicity of \cite{jung-freexentropy}.
\end{proof}

\subsection{$\Delta (F)$ and $\Delta (F:F')$.}

The results of the previous section are insufficient to give a lower
bound for $\Delta (M,\tau )$ and thus for $\beta _{1}^{(2)}(M,\tau )$.
We show, however, that under certain smoothness conditions on the
families $F$ and $F'$, $\Delta (F,F')\geq \Delta (F)$.

\begin{thm}
Let $F=(X_{1},\ldots ,X_{n})$ be a self-adjoint family of generators
of $M$ and let $F'=F\cup (Y_{1},\ldots ,Y_{m})$. Let $D_{1},\ldots ,D_{n}$
be a dual system to $X_{1},\ldots ,X_{n}$ in the sense of \cite{dvv:entropy5};
thus $D_{j}\in B(L^{2}(M))$ satisfy\[
[D_{j},X_{i}]=\delta _{ji}P_{1},\]
where $P_{1}$ denotes the projection onto the trace vector in $L^{2}(M)$.
Assume that $[D_{j},Y_{i}]$ is a Hilbert-Schmidt operator for all
$i$ and $j$. Then\[
\Delta (X_{1},\ldots ,X_{n}:X_{1},\ldots ,X_{n},Y_{1},\ldots ,Y_{m})=\Delta (X_{1},\ldots ,X_{n})=n.\]

\end{thm}
\begin{proof}
Note that for each $j$, the $n+m$-tuple\[
(0,\ldots ,P_{1},\ldots ,0,[Y_{1},D_{j}]^\sigma,\ldots ,[Y_{m},D_{j}]^\sigma)\]
($P_{1}$ in the $j$-th place) lies in $\partial _{F\cup F'}^{\,t}(B(L^{2}(M)))$,
in the notation of (\ref{eq:partialstarF}). Thus\[
\xi _{j}=(0,\ldots ,P_{1},\ldots ,0)\]
($P_{1}$ in the $j$-th place) lies in\[
K=\pi _{n}(\overline{\partial _{F\cup F'}^{\,t}(B(L^{2}(M))}^{w}\cap HS^{n+m}),\]
in the notation of (\ref{eq:dualRelDelta}). Since $(\xi _{1},\ldots ,\xi _{n})$
clearly densely generate $HS^{n}$ as an $M,M$-bimodule, it follows
that the dimension of $K$ over $M\bar{\otimes }M^{o}$ is exactly
$n$. Thus $\Delta (F:F')=n$. Applying this to the case that $m=0$
gives also the estimate for $\Delta (F:F)=\Delta (F)$.
\end{proof}
An important case of existence of a dual system is when $X_{1},\ldots ,X_{n}$
are free semicircular variables; see \cite{dvv:entropy5}.

\section{Appendix I: Abelian von Neumann algebras.}

\n The following theorem is the analog of \cite[Theorem 5.1]{luck:foundations1},
which makes one suspect that its statement should hold more generally
if $A$ is hyperfinite. We were unable to prove this, however. If
the statement holds for $A$ hyperfinite, it would be interesting
if it can be used as a characterization of hyperfinite algebras (see
Remark 5.13 in \cite{luck:foundations1}).

\begin{thm}
\label{pro:TorAvanishes}Let $A$ be a commutative von Neumann algebra
and $\tau$  a normal
faithful trace on $A$. 

 \n (i) Let $f:(A\otimes A^{o})^{n}\rightarrow (A\otimes A^{o})^{m}$ be a
left $A\otimes A^{o}$-module map, then $\beta ^{(2)}(f)=0$.

\n (ii) Let $W$ be an arbitrary $A\otimes A^{o}$-module. Then for all $p\geq 1$,\[
\dim _{A\bar{\otimes }A^{o}}\Tor _{p}^{A\otimes A^{o}}(W,\,A\bar{\otimes }A^{o})=0,\]
\end{thm}

\n Let us first prove a simple lemma,
\begin{lem} \label{sup}
Let $f\in A\otimes A^{o}$, then the spectral projection $p$ of $f^*\,f$ corresponding
to $\ker f$ is the supremum of the projections $e \leq p : e \in A\otimes A^{o}$.
\end{lem} 

\begin{proof}
\n We can assume $A$ is diffuse and identify $A$
with $L^{\infty }([0,1])$, $\tau$ with the Lebesgue measure $\lambda$ 
and $L^{2}(A)$ with $L^{2}([0,1])$. We drop the 
distinction between $A$ and $A^{o}$.

\n Let $f=\sum g_{i}\otimes h_{i}$ and consider $f$ as the function\[
f(x,y)=\sum _{i=1}^{k}g_{i}(x)\,h_{i}(y)\,,\quad x\,,y\in [0,1].\]
Then the projection $p\in A\bar{\otimes }A^{o}=L^{\infty }([0,1]\times [0,1])$
is given by the zero set\[
Z=\{(x,y):f(x,y)=0\}.\]
Recall that a point $z\in Z$ is called a point of density of $Z$ if
the proportion of $Z$ in  squares $S= I\times J$, ($I$ and $J$ intervals
of equal length)
with center $z$ tends to $1$ when their size tends to $0$.
 By Lebesgue's a.e. differentiability theorem, the set $Z$ differs
by a set of measure zero from its set of points of density. We thus
only need to prove the following, with $k$ as above,
\begin{claim}
Let $I,J\subset [0,1]$
be intervals, and $\delta<k^{-1}$ with \[
\frac{\lambda ^{\times 2}((I\times J)\cap Z)}{\lambda ^{\times 2}(I\times J)}>1-\delta^{2}.\]
 Then there are measurable subsets $E\subset I$, $F\subset J$, such
that $E\times F\subset Z$ and\[
\frac{\lambda ^{\times 2}(E\times F)}{\lambda (I\times J)}\geq (1-k\,\delta )^{2}.\]

\end{claim}

\n To prove the claim, let $g:I\rightarrow \mathbb{C}^{k}$,
$h:J\rightarrow \mathbb{C}^{k}$ be given by $
(g(x))_{i}=g_{i}(x)$, $ h(x))_{j}=h_{j}(x)$
 so that,\[
f(x,y)=g(x)\cdot h(y)\,,\quad \forall x\,,y\in [0,1],\]
 where $\cdot $ denotes the standard scalar product on $\mathbb{C}^{k}$.

\n Let
\[
E:= \{x\in I :\lambda \{y:(x,y)\in Z\}>(1-\delta) \lambda (J)\}.\]
Then by Fubini's Theorem $\lambda (E)>(1-\delta) \lambda (I)$.
 Denote by $V(x)$ the subspace of $\mathbb{C}^{k}$ spanned by $g(x)$.
Let $V=\Span (V(x):x\in E)$. Since the dimension of $V$ is at most
$k$, we can choose $x_{1},\ldots ,x_{l}\in E$, $l\leq k$, so that
$V=\Span (g(x_{1}),\ldots ,g(x_{l}))$. For each $1\leq j\leq l$,
the set  $F_{j}$ of $y\in J$ for which $h(y)$ is perpendicular
to $g(x_{j})$ (i.e. $(x_{j},y)\in Z$) has measure at least $(1-\delta )\lambda (J)$.
Thus the measure of $F=\bigcap F_{j}$ is at least $(1-l\delta )\lambda (J)\geq (1-k\,\delta)\lambda (J)$.
But then for all $y\in F$ and $x\in E$, $g(x)\in V$ and $h(y)\perp V$,
so that $f(x,y)=0$. It follows that $E\times F\subset Z$.

\end{proof}
\begin{proof}(of Theorem \ref{pro:TorAvanishes}). First the above lemma implies (i) for $n=m=1$.
Indeed  let $f\in A\otimes A^{o}$, and $p$  the spectral projection of $f^*\,f$ corresponding
to $\ker f$.
The  subspace  $(A\otimes A^{o}) \cdot p $ is dense in $\ker f^{(2)}$. Thus 
since $p$ is a strong limit of projections $e_j \in A\otimes A^{o}$, $p\,e_j=e_j\, p=e_j$,
one gets $(A\otimes A^{o})\cdot e_j\subset \ker(f)$
and the required density of $\ker f$ in $\ker f^{(2)}$.

\n Let now $n$ and $m$ be arbitrary, and reduce to $n=m$ by e.g. replacing $f$ with $f^{*}f$.
Let $F(x,y)$ be the matrix with entries $f_{ij}(x,y)$ and,

\[
F{i_{1}\, \, i_{2}\, \, \cdots \, \, i_{k} \choose j_{1}\, \, j_{2}\, \, \cdots \, \, j_{k}}\]
 the $k\times k$ minor of $F$ obtained by keeping the $i_{1},\ldots ,i_{k}$-th
rows and $j_{1},\ldots ,j_{k}$-th columns of $F$.
Let $Z{i_{1}\cdots i_{k} \choose j_{1}\cdots j_{k}}$ be the zero
set of $F{i_{1}\cdots i_{k} \choose j_{1}\cdots j_{k}}$. Since $F{i_{1}\cdots i_{k} \choose j_{1}\cdots j_{k}}$
is a polynomial expression in the entries of $F$, it belongs to $A\otimes A^{o}$.

\n Let $t\in [0,1]^{2}$ be such that the minors $F{1 \choose 1}$, $\ldots $,
$F{1\cdots r \choose 1\cdots r}$ are all non-zero, while the minors $F{1\cdots r+1 \choose 1\cdots r+1}$,
$\ldots $, $F{1\cdots n \choose 1\cdots n}$ are zero.

\n In this case, the  equation $F\xi =0$, $\xi =(\xi _{1},\ldots ,\xi _{n})$
after performing Gaussian elimination has the form\begin{eqnarray*}
a_{11}^{(0)}\xi _{1}+a_{12}^{(0)}\xi _{2}+\cdots +\cdots +a_{1n}^{(0)}\xi _{n} & = & 0\\
a_{22}^{(1)}\xi _{2}+\cdots +\cdots +a_{2n}^{(1)}\xi _{n} & = & 0\\
\cdots  &  & \\
a_{rr}^{(r-1)}\xi _{r}+\cdots +a_{rn}^{(r-1)}\xi _{n} & = & 0
\end{eqnarray*}
 where\[
a_{ik}^{(p)}=\frac{F{1\cdots p\, \, i \choose 1\cdots p\, \, k}}{F{1\cdots p \choose 1\cdots p}}\]
 (see pp.~24--25 in \cite{gantmaher}). Thus a basis for the null
space of $F$ consists of the vectors $\eta ^{(k,r)}$, $k=1,\ldots ,n-r$,
with coordinates\begin{eqnarray*}
\eta _{t}^{(k,r)} & = & 0,\quad t>r,\, \, t\ne k+r\\
\eta _{r}^{(k,r)} & = & -\frac{a_{rk}^{(r-1)}}{a_{rr}^{(r-1)}},\\
\eta _{r-1}^{(k,r)} & = & -\frac{a_{r-1\, \, k}^{(r-2)}+a_{r-1\, \, r}^{(r-2)}\eta _{r}^{(k)}}{a_{r-1\, \, r-1}^{(r-2)}}\\
 & \cdots  & \\
\eta _{1}^{(k,r)} & = & -\frac{a_{1k}^{(0)}+a_{12}^{(0)}\eta _{2}^{(k)}+\cdots +a_{1r}^{(0)}\eta _{r}^{(k)}}{a_{11}^{(0)}}.
\end{eqnarray*}
 If we set $\xi ^{(k,r)}$ to be the product of $\eta ^{(k,r)}$ by
a sufficiently high power of the (nonzero) expression\[
F{1 \choose 1}\cdots F{1\cdots r \choose 1\cdots r},\]
 we get that for each $k$, the  $\xi _{j}^{(k,r)}$
are polynomials in the entries of $F$, and the vectors $\xi^{(k,r)}$ span
the kernel of $F$.

\n The   polynomial expressions $\xi _{j}^{(k,r)}$ in the entries
of $F$ make sense  without any assumptions on $F$. If
 $F{1\cdots r+1 \choose 1\cdots r+1}$, $\ldots $,
$F{1\cdots n \choose 1\cdots n}$ are zero,  the vectors
$\xi ^{(k,r)}$ lie in the kernel of $F$ (although they will no longer
span the kernel unless $F{1 \choose 1}$, $\ldots $, $F{1\cdots r \choose 1\cdots r}$
are all nonzero).

\n By construction  $\xi _{j}^{(k,r)}\in A\otimes A^{o}$ (as polynomial functions in
$F$). Since $\xi ^{(k,r)}\in \ker F(z)$ for all $z$ such that $F{1\cdots r+1 \choose 1\cdots r+1}$,
$\ldots $, $F{1\cdots n \choose 1\cdots n}$ are zero, using Lemma \ref{sup}
 we therefore get that, \[
\zeta ^{(k,r)}(F)=\xi ^{(k,r)}\chi _{\left\{ z:F{1\cdots r+1 \choose 1\cdots r+1}(z)=\cdots =F{1\cdots n \choose 1\cdots n}(z)=0\right\} }\in \overline{\ker f}\]
 (here $\chi $ denotes the characteristic function of the given set).

\n Applying this result to the 
matrix $F^{\sigma ,\sigma '}$ obtained from $F$ by permuting rows
via a permutation $\sigma $ and columns via a permutation $\sigma '$,
we obtain vectors $$\zeta ^{(k,r,\sigma ,\sigma ')}=\sigma ^{-1}(\xi ^{(k,r)}(F^{\sigma ,\sigma '}))
\in \overline{\ker f}\,.$$

\n For each $z$, let $r$ be  the rank of $F(z)$, we can find $\sigma ,\sigma '$
so that the ${1 \choose 1}$, $\ldots $, ${1\cdots r \choose 1\cdots r}$-minors
of $F^{\sigma ,\sigma '}$ are non-zero and the ${1\cdots r+1 \choose 1\cdots r+1}$,
$\ldots $, ${1\cdots n \choose 1\cdots n}$ minors are zero. Thus
$\{\zeta ^{(k,r,\sigma ,\sigma ')}(z):1\leq k\leq n-r\}$ (and hence
$\{\zeta ^{(k,r,\sigma ,\sigma '}(z)\}_{k,r,\sigma ,\sigma '}$ span
the kernel of $F$ at $z$.
It thus follows that $\ker f$ is dense in $\ker f^{(2)}$ which proves (i).\\

\n Finally the proof of (ii) follows verbatim the argument of Luck (Th. 5.1 \cite{luck:foundations1}).

\end{proof}
\begin{cor}
Let $A$ be an abelian von Neumann algebra. Then for all $k\geq 1$,\[
\beta _{k}^{(2)}(A)=0.\]

\end{cor}
\begin{proof}
By definition,\[
\beta _{k}^{(2)}=\dim _{A\bar{\otimes }A^{o}}\Tor _{k}^{A\otimes A^{o}}(A,A\bar{\otimes }A^{o}),\]
which is zero by the main result of this section.
\end{proof}

\section{Appendix II: Dual Systems.}

\n We recall in this appendix the construction of the dual system in the framework
of section 4, and give the details of the proof of (\ref{comp}).

\n Let $M$ be the von Neumann algebra generated by $X_{1},\ldots ,X_{n}$,
and let $\Omega \in L^{2}(M)$ be the trace vector.
We start by
explicitly constructing the standard form of the von Neumann algebra $N$ 
obtained by adjoining the free semicircular variables $S_{1},\ldots ,S_{n}$.\\ 

\n Consider the vector
space $V=L^{2}(M)\otimes L^{2}(M)\oplus \cdots \oplus L^{2}(M)\otimes L^{2}(M)$
$=(L^{2}(M)\otimes L^{2}(M))^n$,
 and let\[
H=L^{2}(M)\oplus V\oplus (V\otimes _{M}V)\oplus (V\otimes _{M}V\otimes _{M}V)\oplus \cdots \]

\n Note that $V\otimes _{M}V\cong (L^{2}(M)\otimes L^{2}(M)\otimes L^{2}(M))^{n^{2}}$.

\n Then $M$ acts on $H$ both on the right and on the left in the obvious
way, acting on the leftmost or rightmost tensor copy of $V$ each
time. Denote by $\phi _{i}$ the inclusion map from $L^{2}(M)\otimes L^{2}(M)$
into $V$, which places $L^{2}(M)\otimes L^{2}(M)$ as the $i$-th
direct summand.

\n Let $\omega _{i}$ be the $i$-th copy of $1\otimes 1$ in $V$. Denote
by $L_{i}$ and $R_{i}$ the following operators on $H$:\begin{eqnarray*}
L_{i}m & = & \phi _{i}(1\otimes m)\,,\quad \forall m \in L^{2}(M)\\
L_{i}v_{1}\otimes \cdots \otimes v_{n} & = & \omega _{i}\otimes v_{1}\otimes \cdots \otimes v_{n}
\,,\quad \forall v_j \in V\\
R_{i}m & = & \phi _{i}(m\otimes 1)\,,\quad \forall m \in L^{2}(M)\\
R_{i}v_{1}\otimes \cdots \otimes v_{n} & = & v_{1}\otimes \cdots \otimes v_{n}\otimes \omega _{i}
\,,\quad \forall v_j \in V.
\end{eqnarray*}
 Formally,
these are the left and right  tensor multiplications by $\omega _{i}$.

\n It is not hard to check that if we denote by $\lambda $ the left
action of $M$ on $H$, then we have\begin{equation}
L_{i}^{*}\lambda (m)L_{j}=\delta _{ij}\tau (m),\qquad \forall m\in M.\label{eq:Ls}\end{equation}
Similarly, if we denote by $\rho $ the right action of $M$ on $H$,
then we have\[
R_{i}^{*}\rho (m)R_{j}=\delta _{ij}\tau (m),\qquad \forall m\in M.\]
In particular, $L_{i}^{*}L_{i}=R_{i}^{*}R_{i}=1$, and these operators
have norm one.

\n Furthermore,\[
[R_{i},\lambda (m)]=[R_{i}^{*},\lambda (m)]=0,\qquad \forall m\in M.\]
Consider on $B(H)$ the vector state $\psi =\langle \Omega ,\cdot \,\Omega \rangle $,
where $\Omega \in L^{2}(M)$ is regarded as a vector in $H$.
 By a result from \cite{shlyakht:freeness}
 it follows that
if we let\[
S_{i}=L_{i}+L_{i}^{*},\]
then $S_{1},\ldots ,S_{n}$ are a family of free semicircular variables,
free in $(B(H),\psi )$ from $\lambda (M)$. Furthermore, the von
Neumann algebra generated by $M$ (which we identify with $\lambda (M)$)
and $S_{1},\ldots ,S_{n}$ is in standard form, and the operator $J$
is given by\begin{eqnarray*}
J(v_{1}\otimes \cdots \otimes v_{n}) & = & v_{n}^{s }\otimes \cdots \otimes v_{1}^{s },\\
J|_{L^{2}(M)} & = & J_{M},
\end{eqnarray*}
where $s (x\otimes y)=Jy\otimes Jx$, and $J_{M}$is the Tomita
conjugation associated to $M$.

\n Moreover, one has\[
[R_{i},S_{j}]=\delta _{ij}P_{\Omega },\]
which means that $Q_{j}=R_{j}$ is the desired conjugate system to
$S_{1},\ldots ,S_{n}$ relative to $X_{1},\ldots ,X_{n}$, i.e., it
satisfies:\[
[R_{i},X_{j}]=0,\qquad [R_{i},S_{j}]=\delta _{ij}P_{\Omega }.\]
This way, if we set 
$D_{j}(\varepsilon)=E_{W^{*}(X_{1}(\varepsilon),\ldots ,X_{n}(\varepsilon))}\frac{1}{\sqrt{\varepsilon}}Q_{j}E_{W^{*}(X_{1}(\varepsilon),\ldots ,X_{n}(\varepsilon))},$
we get that 
$$
[D_{j}(\varepsilon),X_{i}(\varepsilon)]=\delta _{ij}\,P_{\Omega }
$$
which thus gives the desired
dual system.

\n Moreover, one sees that $J(\sigma (a)Q_{j}\sigma (b))^{*}J1=0$, since
$\sigma (a)^{*}J1$ lies in $L^{2}(M)\subset H$ and $Q_{j}^{*}L^{2}(M)=0$.
On the other hand,\begin{eqnarray*}
\sigma (a)Q_{j}\sigma (b)\cdot 1 & = & \sigma (a)Q_{j}\cdot b=\sigma (a)\cdot (\phi _{j}(b\otimes 1))\\
 & = & J\lambda (a^{*})J\cdot (\phi _{j}(b\otimes 1))\\
 & = & J\lambda (a^{*})\cdot \phi _{j}(1\otimes b^{*})\\
 & = & J\phi _{j}(a^{*}\otimes b^{*})\\
 & = & \phi _{j}(b\otimes a)\\
 & = & (bS_{j}a).
\end{eqnarray*}
which gives (\ref{comp}).
\bibliographystyle{alpha}

\end{document}